\documentclass[a4paper,9pt]{amsart}

\usepackage{amssymb}
\usepackage{color}
\usepackage{mathtools}
\usepackage{amsfonts}
\usepackage{mathrsfs}

\usepackage[colorlinks=true]{hyperref}
\hypersetup{urlcolor=blue, citecolor=red}

\newcommand{\xii}{{\mathrm{|\xi|}}}
\newcommand{\supp}{{\mathrm{\,supp\,}}}
\newcommand{\loc}{{\mathrm{\,loc\,}}}

\newcommand{\eps}{\varepsilon}
\newcommand{\sinc}{{\mathrm{\,sinc\,}}}
\newcommand{\id}{{\mathrm{\,Id\,}}}
\newcommand{\sign}{{\mathrm{\,sign\,}}}

\newcommand{\PV}{\mathrm{P.V.\,}}

\newtheorem{Thm}{Theorem}

\newtheorem{Prop}{Proposition}[section]
\newtheorem{lemma}{Lemma}[section]
\newtheorem{Cor}{Corollary}[section]
\newtheorem{notation}{Notation}
\newtheorem{Rem}{Remark}[section]
\newtheorem{Example}{Example}[section]

\newcommand{\N}{\mathbb{N}}

\newcommand{\R}{\mathbb{R}}
\renewcommand{\S}{\mathbb{S}}

\newcommand{\<}{\langle}
\renewcommand{\>}{\rangle}
\newcommand{\rank}{{\mathrm{rank\,}}}

\title[$L^p-L^q$ estimates for Fourier multipliers]{Sharp $L^p-L^q$ estimates for evolution \\ equations with damped oscillations}

\author[M. D'Abbicco, M.R. Ebert]{Marcello D'Abbicco, Marcelo Rempel Ebert}

\thanks{Orcid: $0000-0003-1369-1157$ (M. D'Abbicco), 0000-0003-4592-252X (M.R. Ebert)}

\address{Marcello D'Abbicco, Department of Mathematics, University of Bari, Via E. Orabona 4, 70125 BARI - ITALY}

\address{Marcelo Rempel Ebert, Departamento de Computa\c{c}\~ao e Matem\'atica, Universidade de S\~ao Paulo, Av. Bandeirantes 3900, Ribeir\~ao Preto, SP, 14040-901, Brasil \& Universit\`a degli Studi di Bari Aldo Moro}

\email{marcello.dabbicco@uniba.it, ebert@ffclrp.usp.br}


\begin{document}

\center{Published on: Mathematische Annalen, January 2025\\
\url{https://doi.org/10.1007/s00208-025-03092-y}}

\begin{abstract}
In this paper we derive sharp $L^p-L^q$ estimates, $1\leq p\leq q\leq \infty$ (including endpoint estimates as $L^1-L^1$ and $L^1-L^\infty$) for dissipative wave-type equations, under the assumption that the dissipation dampen the oscillations but it does not cancel them. We assume that the phase function $w$ is homogeneous of some degree $\sigma>0$ and that its Hessian matrix has maximal rank, including the critical case $\sigma=1$, while the dissipative term $a(\xi)>0$ may be inhomogeneous. The critical case includes waves with viscoelastic or structural damping, damped double dispersion equations and plate equations with rotational inertia, and so on. We also obtain the analogous results for fractional Schr\"odinger-type equations with a potential.
%
\end{abstract}

\maketitle

\section{Introduction}\label{intro}

In this paper, we study $L^p-L^q$ estimates under the effect of a perturbation $Au$ and, respectively, $Au_t$ acting on the initial value problem for first order and, respectively, second order, $\sigma$-evolution equations, where $\sigma>0$.

Namely, we investigate the initial value problem for Schr\"odinger-type equations
\begin{equation}
\label{eq:CPSch}
\begin{cases}
u_t + Au = \pm iL_wu, \quad t>0,\, x\in \R^n,\\
u(0,x)=u_0(x),
\end{cases}
\end{equation}
and for wave-type equations
\begin{equation}
\label{eq:CPwave}
\begin{cases}
u_{tt} + L_{w^2}u + Au_t =0, \quad t>0,\, x\in \R^n, \\
u(0,x)=0, \\
u_t(0,x)=u_1(x).
\end{cases}
\end{equation}
Here
\begin{equation}\label{eq:L}
L_wu = \mathscr{F}^{-1} (w(\xi)\,\hat u),\quad L_{w^2}u = \mathscr{F}^{-1} (w(\xi)^2\,\hat u),
\end{equation}
 $\mathscr{F}$ denotes the Fourier transform in the space of tempered distributions~$\mathcal S'$ and $\mathscr{F}^{-1}f(x)=(2\pi)^{-n}\mathscr{F}f(-x)$ its inverse.

We assume that  the phase $w$ is real-valued, homogeneous of degree $\sigma>0$ and smooth on $\R^n\setminus\{0\}$ (in most cases, $w\in\mathcal C^{n+3}(\R^n\setminus\{0\})$ shall be enough). We also assume that the Hessian matrix $H_w$ has maximal rank, that is, $\det H_w\neq0$ if $\sigma\neq1$ and $\rank H_w=n-1$ if $\sigma=1$ (since the homogeneity of degree $1$ implies that $\det H_w=0$). In~\eqref{eq:CPwave} we also assume that $w(\xi)>0$ for any $\xi\neq0$. 
We define
\[ Au=\mathscr{F}^{-1}(a(\xi)\hat u), \]
with $a\in\mathcal C^{n+1}(\R^n\setminus\{0\})$, such that $a(\xi)>0$ for $\xi\neq0$, and that $a(\xi)$ verifies a suitable control with its derivatives, at low and high frequencies.

The problem in~\eqref{eq:CPSch} includes (perturbations of) the classical Schr\"odinger equation, when $w=\xii^2$, as well as fractional Schr\"odinger equations or higher order Schr\"odinger equations, when $w=\xii^\sigma$; it also includes the (linear) Korteweg-de Vries equation $u_t+\partial_x^3u=0$, for which $w=-\xi^3$. The problem in~\eqref{eq:CPwave} includes (perturbations of) the wave equation when $w=\xii$, the plate equation when $w=\xii^2$, the double dispersion equation or plate equation with rotational inertia, when $A=(\id-\Delta)^{-1}B$ for some dissipative operator $B$, and many other models. Equations as the ones in~\eqref{eq:CPSch} and~\eqref{eq:CPwave} are called  $\sigma$-evolution equations, respectively of first order and second order (in time). This name goes back to~\cite{Miz61}; several results of well-posedness for $\sigma$-evolution differential equations have been obtained in recent years, see for instance, \cite{AAC2024,AB2013,ABZ2012,CHR2008}.

When $A\equiv0$, the study of $L^p-L^q$ estimates for \eqref{eq:CPSch} and~\eqref{eq:CPwave} mostly reduces to the study of the multipliers $(1-\chi)\,e^{\pm iw(\xi)}$, where $\chi\in\mathcal C_c^\infty$ verifies $\chi=1$ in a neighborhood of the origin. Since $w(\xi)$ is homogeneous, the decay estimates depend on the homogeneity degree $\sigma$ of $w$ (see later, Theorems~\ref{thm:A0} and~\ref{thm:A0w}). This study goes back to~\cite{Miy81,P,Sjo} and to many other authors.

The assumption that $\rank H_w=n-1$ when $\sigma=1$ or $\rank H_w=n$ when $\sigma\neq1$, deeply influences dispersive estimates of type $L^p-L^{p'}$, with $p'=p/(p-1)$; roughly speaking, dispersive estimates are worse when $H_w$ is more singular. On the other hand, the regularity theory, i.e., $L^p-L^p$ estimates, is better when $\sigma=1$ (for $w=\xii$, see~\cite{Sjo,Wai}). Therefore, the case $\sigma=1$ can be considered critical, with respect to $\sigma\neq1$.

In both the problems~\eqref{eq:CPSch} and~\eqref{eq:CPwave}, if the perturbation acts at low frequencies, it influences the decay estimates, whereas if it acts at high frequencies, it influences the regularity of the solution. This influence is limited to some region of the $(p,q)$-plane (see Remarks~\ref{Rem:noASch} and~\ref{Rem:noAw}). In~\eqref{eq:CPSch}, the influence of the perturbation $Au$ is related to the interplay of the two multipliers, $e^{\pm iw(\xi)t}$ and $e^{-a(\xi)t}$. In~\eqref{eq:CPwave}, the action of the dissipation is more complicated. In the damped oscillations regime (see later, \eqref{eq:wavesol}), the dissipation does not cancel the oscillations and it modifies the phase function, so that the study of the interplay of the multipliers is more complicated. Therefore, we use a perturbation argument instead of a direct study of the multiplier. 

Obtaining $L^p-L^q$ estimates in the full range $1\leq p\leq q\leq \infty$, for Cauchy problems for evolution equations, provides a useful tool to study problems under the effect of a nonlinear perturbation. In particular, the presence of a damping term in the damped oscillations regime hints to the possibility to mix decay estimates of diffusive type, in particular $L^1-L^q$ estimates, and of dispersive type, $L^p-L^{p'}$. The first kind of estimates are more related to diffusive problems and generate Fujita type critical exponents (see, for instance, \cite{M76,TY01}, see also~\cite{DAR14,EGR2020,PKR15}), but this type of critical exponent also appears for model as~\eqref{eq:CPwave} for $\sigma\neq1$ in low dimension, see~\cite{DAE22}. 
The second type of estimates are more related to hyperbolic problems and their $\sigma$-evolution counterparts, and often involve the so-called Strauss exponent~\cite{Strauss,S89}, that appears in the wave equation (see, for instance, \cite{Georgiev,GeorgLinbSogge1997,Glassey1981Ex,John1979,LS96,Schaeffer1985,Zhou1995}, see also~\cite{CR24}), in the Boussinesq equation~\cite{CO07}, and in several other models. We expect that our $L^p-L^q$ estimates in this paper may lead to some critical exponents that is someway intermediate between the two scenarios.

\subsection{Estimates for Fourier multipliers and homogeneous evolution equations}

Following~\cite{Ho}, in this paper we say that a multiplier $m$ is in $M_p^q$ if $m\in\mathcal S'(\R^n)$ (space of tempered distributions), and
\begin{equation}\label{eq:Mpq}
\|m\|_{M_p^q}=\sup\big\{\|\mathscr{F}^{-1}(m\mathscr{F}(f))\|_{L^q}:f\in \mathcal{S}(\R^n), \|f\|_{L^p}=1\big\},
\end{equation}
is finite. Here $\mathscr{F}$ denotes the Fourier transform with respect to $x$ and we write $\hat f=\mathscr{F}f$ for any $f\in \mathcal{S}'$. %
The duality property $M_p^q=M_{q'}^{p'}$ holds, where $p'=p/(p-1)$ denotes the H\"older conjugate exponent of $p$. 
%
%
We also define $M(H^1,L^q)$ and $M(H^1,H^1)$ as the space of bounded multipliers from the real Hardy space $H^1$ to $L^q$ and to itself. It is convenient for us to describe the real Hardy space $H^1$ as the subset of $L^1$ functions~$f$ such that the Riesz transforms of $f$ are also in $L^1$; we put
\begin{equation}\label{eq:H1norm}
\|f\|_{H^1} = \|f\|_{L^1} + \sum_{j=1}^n \|R_jf\|_{L^1},
\end{equation}
where $R_j$ are the Riesz transforms of $f$ defined by $R_jf=\mathscr{F}^{-1}(i(\xi_j/\xii)\,\hat f)$ (see~\cite{F71}). Therefore,
\begin{equation}\label{eq:H1multnorms}\begin{split}
\|m\|_{M(H^1,L^q)}&=\sup\big\{\|\mathscr{F}^{-1}(m\mathscr{F}(f))\|_{L^q}:f\in \mathcal{S}\cap H^1, \|f\|_{H^1}=1\big\}, \\
\|m\|_{M(H^1,H^1)}&=\sup\big\{\|\mathscr{F}^{-1}(m\mathscr{F}(f))\|_{H^1}:f\in \mathcal{S}\cap H^1, \|f\|_{H^1}=1\big\}.
\end{split}\end{equation}
%
%
%
For any $1\leq p\leq q\leq\infty$, we define
\begin{equation}\label{eq:dpqdef}\begin{split}
d(p,q)
    & = \frac{n}\sigma\left(\frac1p-\frac1q\right)+n\max\left\{\frac12-\frac1p,\frac1q-\frac12\right\} \qquad \sigma\neq1,\\
d(p,q)
    & = n\left(\frac1p-\frac1q\right)+(n-1)\max\left\{\frac12-\frac1p,\frac1q-\frac12\right\} \qquad \sigma=1.
\end{split}\end{equation}
%
The dual property $d(q',p')=d(p,q)$ is related to the dual property of $M_p^q$. Then we have the following.

\begin{Thm}\label{thm:w}
Let $w$ be homogeneous of degree $\sigma$ and with maximal rank, and assume that $\phi\in S^{-b\sigma}$, that is, $\phi\in\mathcal C^\infty$ and
\[ |\partial_\xi^\alpha \phi(\xi)|\leq C_\alpha\,\<\xi\>^{-b\sigma-|\alpha|}. \]
If $b\geq d(p,q)$, then
\begin{equation}\label{eq:wbasic}\begin{split}
& \phi \,e^{\pm iw(\xi)} \in M_p^q, \quad 1<p\leq q<\infty,\\
& \phi \,e^{\pm iw(\xi)} \in M(H^1,L^q), \quad q\in[1,\infty].
\end{split} \end{equation}
Moreover, if $\sigma\neq1$ and $b\geq d(1,\infty)=n(2-\sigma)/2\sigma$, then
\begin{equation}\label{eq:endpoint}
\phi \,e^{\pm iw(\xi)} \in M_1^\infty.
\end{equation}
If $\phi$ vanishes in a neighborhood of $\xi=0$ and $b>d(p,q)$, then
\begin{equation}\label{eq:w1}
\phi\,e^{\pm iw(\xi)} \in M_1^q, \quad q\in[1,\infty].
\end{equation}
Without the maximal rank assumption, the previous estimates hold for $1\leq p\leq q\leq2$ and for its dual range $2\leq p\leq q\leq\infty$.
\end{Thm}
The proof of Theorem~\ref{thm:w} for $\sigma=1$ is contained in~\cite{Stein93}, see chapter IX, \textsection 6.16 at p. 428, and it extends the result in~\cite{P} for $w=\xii$. Theorem~\ref{thm:w} is proved when $\sigma\neq1$ for $w=\xii^\sigma$ in~\cite{Miy81}, we prove it in the general case in \textsection\ref{sec:w}. The condition on $b$ is sharp, in general (see~\cite{Miy81,P} for $w=\xii^\sigma$).

Theorem~\ref{thm:w} leads to the following immediate consequences, for the solutions to~\eqref{eq:CPSch} and~\eqref{eq:CPwave}, when $A=0$.
\begin{Thm}\label{thm:A0}
Let $1\leq p\leq q\leq\infty$ and $d(p,q)$ as in~\eqref{eq:dpqdef}. Let $s=(d(p,q))_+$ if $1<p\leq q<\infty$ or let $s\geq0$ be such that $s>d(p,q)$ otherwise. If $u_0\in H^{s\sigma,p}(\R^n)$, then $\hat u(t,\xi) = e^{\pm itw(\xi)}\,\hat u_0(\xi)$ verifies the estimate
\begin{equation}\label{eq:E1}
\|u(t,\cdot)\|_{L^q} \leq C\,t^{-\frac{n}\sigma\left(\frac1p-\frac1q\right)}\,(1+t)^s\,\|u_0\|_{H^{s\sigma,p}},\qquad t>0.
\end{equation}
Estimate~\eqref{eq:E1} also holds with $s=(d(1,\infty))_+$ if $(p,q)=(1,\infty)$ and $\sigma\neq1$.
\end{Thm}
For the optimality of estimate~\eqref{eq:E1} when $w=c\xii^\sigma$, see, for instance, \cite[Propositions 7.1, 7.3]{DAE21}.

In the radial case $w(\xi)=w(|\xi|)$, \eqref{eq:CPSch} with $A\equiv0$ is studied in~\cite{Ozawa2011} for $w$ not necessarily homogeneous. For more about dispersive estimates for fractional Schr\"odinger equations, see also~\cite{Guo2008, Laskin} and the references therein. As far as we know, result analogous to the one in Theorem~\ref{thm:A0} are not widely studied when $w$ is non radial. Even in the case of radial homogeneous $w$, that is, $w=c\xii^\sigma$, the ``critical'' case $\sigma=1$ is excluded in~\cite{Ozawa2011}.

The solution to~\eqref{eq:CPwave} in appropriate spaces verifies
\begin{equation}\label{eq:uwave}
\hat u(t,\xi) = t\,\sinc (tw(\xi))\,\hat u_1(\xi),
\end{equation}
where $\sinc \rho = \rho^{-1}\sin\rho$ is the cardinal sine function.
\begin{Thm}\label{thm:A0w}
Let $1\leq p\leq q\leq\infty$ and $d(p,q)$ as in~\eqref{eq:dpqdef}. Let $s=(d(p,q)-1)_+$ if $1<p\leq q<\infty$ or let $s\geq0$ be such that $s>d(p,q)-1$ otherwise. If $u_1\in H^{s\sigma,p}(\R^n)$, then~\eqref{eq:uwave} verifies the estimate
\begin{equation}\label{eq:A0w}
\|u(t,\cdot)\|_{L^q} \leq C\,t^{1-\frac{n}\sigma\left(\frac1p-\frac1q\right)}\,(1+t)^s\,\|u_1\|_{H^{s\sigma,p}},\qquad t>0.
\end{equation}
Estimate~\eqref{eq:A0w} also holds with $s=(d(1,\infty)-1)_+$ if $(p,q)=(1,\infty)$ and $\sigma\neq1$. Moreover, if $w=c\xii$, one may take $s=0$ in the following two exceptional cases: if $n=1$ and $(p,q)=(1,\infty)$, and if $n=3$ and $p=q=1$ or $p=q=\infty$.
\end{Thm}
For the optimality of estimate~\eqref{eq:A0w} when $w=c\xii^\sigma$, see, for instance, \cite[Propositions 7.1, 7.3]{DAE21}. The regularity assumption is the same obtained for more general hyperbolic equation in~\cite[Theorem 4.1]{SSS} for the case $q=p\in (1,\infty)$.

\subsection{Results for perturbed fractional Schr\"odinger equations}\label{sec:Schr}

We consider the problem for the perturbed fractional Schr\"odinger equation~\eqref{eq:CPSch}, with $L$ as in~\eqref{eq:L}. The solution to~\eqref{eq:CPSch} in appropriate spaces verifies
\begin{equation}\label{eq:uSch}
\hat u(t,\xi) = e^{\pm itw(\xi)-ta(\xi)}\,\hat u_0(\xi).
\end{equation}
Our aim is to show that when $d(p,q)>0$, the presence of the potential $Au$ may improve the previous estimate in different ways. The effect of $A$ is different at low and high frequencies, so we distinguish two cases. To localize functions at low or high frequencies,  we use the notation
\[  T_\chi f=\mathscr{F}^{-1}(\chi\hat f), \]
where $\chi\in\mathcal C_c^\infty$ is a cut-off function with $\chi=1$ in a neighborhood of the origin (see later, Notation~\ref{not:chi}).

We assume that $a\in \mathcal C^{n+1}(\R^n\setminus\{0\})$ and that $a(\xi)>0$ for any $\xi\neq0$. When we study a problem at low frequencies, we assume that
\begin{equation}
\label{eq:oscillationslow}
a(\xi)\geq a_0\,\xii^{\theta_0},\quad |\partial_\xi^\gamma a(\xi)|\leq C\,\xii^{\theta_0-|\gamma|},\quad |\gamma|\leq n+1,
\end{equation}
for some given $\theta_0>0$, $a_0>0$ and $C>0$, in a neighborhood of the origin. When we study a problem at high frequencies, we assume that
\begin{equation}
\label{eq:oscillationshigh}
a(\xi)\geq a_1\,\xii^{\theta_1},\quad |\partial_\xi^\gamma a(\xi)|\leq C\,\xii^{\theta_1-|\gamma|},\quad|\gamma|\leq n+1,
\end{equation}
for some given $\theta_1>0$, $a_1>0$ and $C>0$, out of some compact set of $\R^n$. 
%
%
In particular, if $a\in\mathcal C^{n+1}(\R^n\setminus\{0\})$ is homogeneous of degree~$\theta>0$, \eqref{eq:oscillationslow} and~\eqref{eq:oscillationshigh} follow with $\theta_0=\theta_1=\theta$, as a consequence of the homogeneity.

At low frequencies, the potential $Au$ improves the decay rate with respect to the one in~\eqref{eq:E1}. Clearly, the regularity does not come into play, since the solution is localized at low frequencies.
\begin{Thm}\label{thm:Schdec}
Let $1\leq p\leq q\leq\infty$ be such that $d=d(p,q)>0$,where $d(p,q)$ is as in~\eqref{eq:dpqdef}. Moreover, let us assume \eqref{eq:oscillationslow} and  that $u_0\in L^p(\R^n)$. If $\theta_0>\sigma$, then~\eqref{eq:uSch} verifies the following decay estimate:
\begin{equation}\label{eq:Schdec}
\|T_\chi\,u(t,\cdot)\|_{L^q} \leq C\,(1+t)^{-\frac{n}\sigma\left(\frac1p-\frac1q\right)+d\left(1-\frac\sigma{\theta_0}\right)}\,\|u_0\|_{L^p}, \quad t\geq0.
\end{equation}
If $\theta_0\in(0,\sigma]$, then we have the following decay estimate:
\[
\|T_\chi\,u(t,\cdot)\|_{L^q} \leq C\,(1+t)^{-\frac{n}\sigma\left(\frac1p-\frac1q\right)}\,\|u_0\|_{L^p}, \quad t\geq0.
\]
%
\end{Thm}
Comparing~\eqref{eq:Schdec} with~\eqref{eq:E1}, the term $(1+t)^{-d\frac\sigma{\theta_0}}$ corresponds to the extra decay rate associated to the potential $Au$ when $\theta_0>\sigma$. When $\theta_0\leq\sigma$, the extra decay rate is sufficient to cancel the loss of decay $(1+t)^s$ associated to oscillations in Theorem~\ref{thm:A0}.

%
%
%
At high frequencies, the potential $Au$ allows to reduce the regularity $s\sigma$ that appears in Theorem~\ref{thm:A0} for the initial data, thought a additional singular power appears as $t\to0$ if $\theta_1<\sigma$, with respect to $t^{-\frac{n}\sigma\left(\frac1p-\frac1q\right)}$.
\begin{Thm}\label{thm:Schreg}
Let $1\leq p\leq q\leq\infty$ b be such that $d=d(p,q)>0$, where $d(p,q)$ is as in~\eqref{eq:dpqdef}. Moreover, let us assume \eqref{eq:oscillationshigh}. If~$\theta_1\in(0,\sigma)$, assume that $u_0\in H^{s\sigma,p}(\R^n)$, where
\[ 0 \leq s \leq \left(1-\frac{\theta_1}\sigma\right)\,d.\]
Then~\eqref{eq:uSch} verifies the following estimate:
\begin{equation}\label{eq:Schreg}
\|(\id-T_\chi)\,u(t,\cdot)\|_{L^q} \leq C\,t^{-\frac{n}\sigma\left(\frac1p-\frac1q\right)+d-(d-s)\frac\sigma{\theta_1}} \,e^{-ct}\,\|u_0\|_{H^{s\sigma,p}},\quad t>0.
\end{equation}
If $\theta_1\geq\sigma$, assume that $u_0\in L^p(\R^n)$. Then,
%
\[ \|(\id-T_\chi)\,u(t,\cdot)\|_{L^q} \leq C\,t^{-\frac{n}\sigma\left(\frac1p-\frac1q\right)} \,e^{-ct}\,\|u_0\|_{L^p},\quad t>0.\]
%
\end{Thm}
In particular, when $q=p$, the estimate in Theorem~\ref{thm:Schreg} is not singular at $t=0$ if $\theta_1\geq\sigma$ or if
\[ \theta_1\in(0,\sigma), \qquad d(p,p) \leq \frac{s}{1-\frac{\theta_1}\sigma}. \]
\begin{Rem}\label{Rem:noASch}
When $d(p,q)\leq 0$, the potential $Au$ has no influence in the $L^p-L^q$ estimates, namely, $T_\chi\,u$ verifies~\eqref{eq:E1} for large $t$ and $(\id-T_\chi)\,u$ verifies~\eqref{eq:E1} for small $t$, with $s$ as in Theorems~\ref{thm:Schdec} and~\ref{thm:Schreg}. This means that $s=0$, with the exception of $s>0$ in the limit case $d(p,q)=0$ with $p=1$ or $q=\infty$ (exception given for the special case $(p,q)=(1,\infty)$, and $\sigma\neq1$, since we may take $s=0$ in that case).
\end{Rem}
The estimates in Theorems~\ref{thm:Schdec} and~\ref{thm:Schreg} are sharp, in the sense discussed in~\textsection\ref{sec:optimality}.

%
Theorems~\ref{thm:Schdec} and~\ref{thm:Schreg} may be easily combined if $a(\xi)$ is homogeneous of degree $\theta$. Let $d(p,q)>0$. Then:
\[ \|u(t,\cdot)\|_{L^q} \leq \begin{cases}
C\,t^{-\frac{n}\sigma\left(\frac1p-\frac1q\right)}\,(1+t)^{d\left(1-\frac\sigma{\theta}\right)}\,\|u_0\|_{L^p}, & \text{if $\theta\geq\sigma$,}\\
C\,t^{-\frac{n}\sigma\left(\frac1p-\frac1q\right)}\,\big(1+t^{-\left((d-s)\frac\sigma{\theta}-d\right)_+}\big)\,\|u_0\|_{H^{s\sigma,p}},& \text{if $\theta\in(0,\sigma]$.}
\end{cases} \]

\subsection{Results for perturbed wave equations}\label{sec:wave}

The study of the effect of the perturbation $Au$ on the Schr\"odinger equation~\eqref{eq:CPSch} is relatively easy using the theory of multipliers. On the other hand, the study of the initial value problem for the wave equation~\eqref{eq:CPwave} is more complicated, especially when $\sigma=1$. This is due to the fact that the phase function in the fundamental solution is modified by the presence of the dissipation $A$, see later, \eqref{eq:wavesol}.

To manage this perturbation in the phase, 
we expand the fundamental solution to~\eqref{eq:CPwave} by Taylor's theorem in order to isolate the main term with the unperturbed phase function $w$. The remaining terms are managed by simpler multiplier lemmas.

As in \textsection\ref{sec:Schr}, we are interested in studying how the dissipation $Au_t$ influences estimates at low and at high frequencies. However, now the situation is more complicated. Indeed, performing the Fourier transform on the equation in~\eqref{eq:CPwave}, we obtain the ODE, depending on the parameter $\xi$:
\[ \hat u_{tt} + w(\xi)^2 \hat u + a(\xi)\hat u_t=0, \]
that is, the ODE of the damped harmonic oscillator, with initial data $\hat u(0,\xi)=0$ and $\hat u_t(0,\xi)=\hat u_1$. 
%
%
There are two possible regimes, according to the sign of $a(\xi)^2-4w(\xi)^2$:
\begin{itemize}
    \item when $a(\xi)^2<4w(\xi)^2$, we are in the ``damped oscillations regime'': the attrition term $a(\xi)$ dampen the oscillations, without canceling them, and the solution is
    \begin{equation}\label{eq:wavesol}
    \hat u (t,\xi) = e^{-a(\xi)\frac{t}2}\,\sinc (w(\xi)\sqrt{1-\tilde a(\xi)}\,t)\,\hat u_1, \quad \text{where}\qquad \tilde a(\xi)=\frac{a(\xi)^2}{4w(\xi)^2};
    \end{equation}
    \item when $a(\xi)^2\geq4w(\xi)^2$, we are in the ``overdamping regime'': the attrition term $a(\xi)$ cancels the oscillations, and the solution is
    \begin{equation}\label{eq:wavesolOD}
    \hat u (t,\xi) =\frac{e^{\lambda_+t}-e^{\lambda_-t}}{\lambda_+-\lambda_-}\,\hat u_1,\quad \text{where}\qquad \lambda_\pm = \frac{-a(\xi)\pm \sqrt{a(\xi)^2-4w(\xi)^2}}2\,,
    \end{equation}
    where~\eqref{eq:wavesolOD} is intended as its limit $te^{-w(\xi)t}\hat u_1$ when $a(\xi)=2w(\xi)$.
\end{itemize}
We are interested in the ``damped oscillations regime'', since in the ``overdamping regime'' oscillations are canceled and the multipliers are exponentials with negative real part, i.e., of diffusive type.

Since we are only interested in frequencies when $\xi\to0$ and when $\xii\to\infty$, in order to be in the ``damped oscillations regime'', it is sufficient to assume that $\theta_0>\sigma$ to get that $a(\xi)^2<4w(\xi)^2$ in some neighborhood of the origin, and that $\theta_1<\sigma$ to get that $a(\xi)^2<4w(\xi)^2$ out of some compact set. We stress that the assumption that $w(\xi)>0$ never vanishes for $\xi\neq0$ and that it is homogeneous, is crucial here.

We are now ready to state our results. As we did in \textsection\ref{sec:Schr}, we distinguish low and high frequencies, and we provide analogous results to Theorems~\ref{thm:Schdec} and~\ref{thm:Schreg}.
\begin{Thm}\label{thm:Wavedec}
Let $1\leq p\leq q\leq\infty$ be such that $d=d(p,q)>1$, where $d(p,q)$ is as in~\eqref{eq:dpqdef}. Moreover, let us assume \eqref{eq:oscillationslow} with $\theta_0>\sigma$, and  that $u_0\in L^p(\R^n)$. Then~\eqref{eq:wavesol}-\eqref{eq:wavesolOD} verifies the following decay estimate:
\begin{equation}\label{eq:Wavedec}
\|T_\chi\,u(t,\cdot)\|_{L^q} \leq C\,(1+t)^{1-\frac{n}\sigma\left(\frac1p-\frac1q\right)+(d-1)\left(1-\frac\sigma{\theta_0}\right)}\,\|u_1\|_{L^p}, \quad t\geq0.
\end{equation}
%
\end{Thm}
The term $(1+t)^{-(d-1)\frac\sigma{\theta_0}}$ corresponds to the extra decay rate associated to the dissipation $Au_t$, with respect to~\eqref{eq:A0w}.
\begin{Thm}\label{thm:Wavereg}
Let $1\leq p\leq q\leq\infty$ b be such that $d=d(p,q)>1$, where $d(p,q)$ is as in~\eqref{eq:dpqdef}. Moreover, let us assume \eqref{eq:oscillationshigh} with~$\theta_1\in(0,\sigma)$, and that $u_0\in H^{s\sigma,p}(\R^n)$, where
\[ 0 \leq s \leq \left(1-\frac{\theta_1}\sigma\right)\,(d-1).\]
Then~\eqref{eq:wavesol}-\eqref{eq:wavesolOD} verifies the following estimate:
\begin{equation}\label{eq:Wavereg}
\|(\id-T_\chi)\,u(t,\cdot)\|_{L^q} \leq C\,t^{1-\frac{n}\sigma\left(\frac1p-\frac1q\right)+(d-1)-(d-1-s)\frac\sigma{\theta_1}} \,e^{-ct}\,\|u_1\|_{H^{s\sigma,p}},\quad t>0.
\end{equation}
%
\end{Thm}
In particular, when $q=p$, the estimate in Theorem~\ref{thm:Wavereg} is not singular at $t=0$ if
\[ 1<d(p,p) \leq 1+\frac{s}{1-\frac{\theta_1}\sigma}. \]
%
%
The estimates in Theorems~\ref{thm:Wavedec} and~\ref{thm:Wavereg} are sharp, in the sense discussed in~\textsection\ref{sec:optimality}.
\begin{Rem}\label{Rem:noAw}
When $d(p,q)\leq 1$, the dissipation $Au_t$ has no influence in the $L^p-L^q$ estimates, namely, $T_\chi\,u$ verifies~\eqref{eq:A0w} for large $t$ and $(\id-T_\chi)\,u$ verifies~\eqref{eq:A0w} for small $t$, with $s$ as in Theorems~\ref{thm:Wavedec} and~\ref{thm:Wavereg}. This means that $s=0$, with the exception of $s>0$ in the limit case $d(p,q)=1$ with $p=1$ or $q=\infty$ (exception given for the special case $(p,q)=(1,\infty)$, and $\sigma\neq1$, since we may take $s=0$ in that case).
\end{Rem}

\subsubsection{The $L^2$ theory}\label{sec:L2}

When $p=2$ or $q=2$, oscillations become irrelevant in the multipliers, and a much simpler proof can be provided for Theorems~\ref{thm:Schdec} and~\ref{thm:Schreg} as well as for Theorems~\ref{thm:Wavedec} and~\ref{thm:Wavereg}, directly using Plancherel theorem or Haussdorff-Young inequality. Moreover, for \eqref{eq:CPwave}, when $(p,q)=(1,2)$ or, by duality, $(p,q)=(2,\infty)$, this may provide a benefit in space dimension $n=2\sigma$, with respect to the estimates provided by Theorem~\ref{thm:Wavedec} and by Theorem~\ref{thm:Wavereg} when~$s=0$. Indeed, $d(1,2)=n/(2\sigma)$ so Theorem~\ref{thm:A0w} produces a loss of decay $(1+t)^\varepsilon$, where $\varepsilon=s-n/(2\sigma)>0$. This arbitrarily small polynomial loss may be relaxed to a logarithmic loss of decay. Similarly, the regularity $H^{\varepsilon\sigma,1}$ may be reduced to $L^1$, paying a logarithmic singular power as $t\to0$.
\begin{Prop}\label{prop:Wavedec2}
Let $n=2\sigma$, $\theta_0>\sigma$, and assume that $u_1\in L^1$. Then~\eqref{eq:wavesol} verifies the following estimate:
\begin{equation}\label{eq:Wavedec2}
\|T_\chi\,u(t,\cdot)\|_{L^2} \leq C\,(1+(\log t)_+)^{\frac12}\,\|u_1\|_{L^1}, \quad t\geq0.
\end{equation}
\end{Prop}
\begin{Prop}\label{prop:Wavereg2}
Let $n=2\sigma$, $\theta_1\in(0,\sigma)$, and assume that $u_1\in L^p$. Then~\eqref{eq:wavesol} verifies the following estimate:
\begin{equation}\label{eq:Wavereg2}
\|(\id-T_\chi)\,u(t,\cdot)\|_{L^2} \leq C\,(1+(-\log t)_+)^{\frac12}\,e^{-ct}\,\|u_1\|_{L^1},\quad t>0.
\end{equation}
\end{Prop}
The analogous results of Propositions~\ref{prop:Wavedec2} and~\ref{prop:Wavereg2} hold for $(p,q)=(2,\infty)$.

The estimates in Proposition~\ref{prop:Wavedec2} and~\ref{prop:Wavereg2} are sharp, see later, Lemma~\ref{lem:crucialwave2}.
\begin{Rem}\label{rem:wave}
When $w=c\xii$, we have
\[ d(1,q)=n\left(1-\frac1q\right)+(n-1)\left(\frac1q-\frac12\right), \]
so the equality $d(1,q)=1$ mentioned in Remark~\ref{Rem:noAw} corresponds to the three special cases; the case $n=1$ and $(p,q)=(1,\infty)$, and the case $n=3$ and $p=q=1$, were already included in Theorem~\ref{thm:A0w}. The third special case is $n=2$ and $(p,q)=(1,2)$, and it is included in Propositions~\ref{prop:Wavedec2} and~\ref{prop:Wavereg2}.

This gives a complete picture of sharp estimates for the wave model, while for other models the optimality in the limit case $d(p,q)=1$ and $p=1$ or $q=\infty$ remains open (expect for $(p,q)=(1,\infty)$ when $\sigma\neq1$).
\end{Rem}
We stress that in this paper we may consider very general dissipative operators, as the following example shows.
\begin{Example}
Let $T$ be the sum of a positive constant $c_0$ with a singular integral operator:
\[ Tf = c_0f + \PV \frac{\Omega(x/|x|)}{|x|^n}\ast f, \]
where $\Omega\in L^1(\S^{n-1})$, has zero average. 
Then its symbol
\[ m(\xi)=c_0 + \int_{\S^{n-1}} \Gamma((\xi\cdot y)/\xii)\,\Omega(y)\,dS(y), \quad \Gamma(t)= -\frac{i\pi}2\,\sign t - \log|t|, \]
is homogeneous of degree zero. If, moreover, $\Omega\in\mathcal C^\infty(\S^{n-1})$, then $m\in\mathcal C^\infty(\R^{n}\setminus\{0\})$ (see, for instance, \cite[chapter 3, \textsection3.5]{Stein70}) and the converse also holds. 
Suitable linear combinations of Riesz transforms and higher order Riesz transforms fit in this scenario.

In particular, if $m(\xi)>0$, taking the composition of $T$ with operators of type $(\id-\Delta)^\eta(-\Delta)^{\frac{\theta}2}$ with $\theta>0$ and $\eta>-\theta/2$ generates examples of $A=(\id-\Delta)^\eta(-\Delta)^{\frac{\theta}2}T$ with $\theta_0=\theta$ in~\eqref{eq:oscillationslow} and $\theta_1=\theta+2\eta$ in~\eqref{eq:oscillationshigh}.

We stress that $A$ is anisotropic, in general, for instance, we may consider
\[ A = \sum_{jk} c_{jk} \partial_{x_j}\partial_{x_k}, \]
provided that the symbol is elliptic, i.e., $a(\xi)=\sum_{jk} c_{jk} \xi_j\xi_k >0$.

We may also add perturbations to $a(\xi)$ which possibly contain fast oscillations with respect to~$\xi$. For instance, if $a(\xi)$ verifies \eqref{eq:oscillationshigh} for some $\theta_1>0$, then
\[ a(\xi) + \sin (\xii^{-\eta}),  \]
also verifies \eqref{eq:oscillationshigh}, provided that $(n+1)(1-\eta)\leq\theta_1$.
\end{Example}

\subsection{Plan of the paper}

The plan of the paper is the following:
\begin{itemize}
\item in \textsection\ref{sec:examples}, we discuss the application of our results to several dissipative wave models;
\item in \textsection\ref{sec:multipliers} we recall some known results on the theory of $M_p^q$ multipliers and we state several lemmas which were used to deal with the multiplier $e^{-a(\xi)t}$ in the paper (postponing their proofs to~\textsection\ref{sec:multproofs});
\item in \textsection\ref{sec:w} we prove Theorem~\ref{thm:w} and its direct consequences: Theorems~\ref{thm:A0} for the fractional Schr\"odinger equation, and Theorem~\ref{thm:A0w} for the wave-type equation;
\item in \textsection\ref{sec:Sch}, we prove Theorems~\ref{thm:Schdec} and~\ref{thm:Schreg} for the fractional Schr\"odinger equation with potential;
\item in \textsection\ref{sec:waveproofs}, we prove Theorems~\ref{thm:Wavedec} and~\ref{thm:Wavereg} for the wave-type equation with dissipation;
\item in \textsection\ref{sec:optimality}, we prove the optimality of our $L^p-L^q$ estimates;
\item in \textsection\ref{sec:multproofs}, we provide the proofs of multiplier Lemmas introduced in \textsection\ref{sec:multipliers} to deal with $e^{-a(\xi)t}$.
\end{itemize}

\subsection{Notation}

In this paper, we use the following notation.

\begin{notation}\label{lesssim}
Let $f,g : [0,\infty)\to(0,\infty)$ be two  positive functions. We use the notation $f\lesssim g$ if there exists a positive constant  $C$ such that $  f(t) \leq C g(t)$, for all $t\geq0$.
\end{notation}
\begin{notation}
We use the notation $\sin^2(\xii)=(\sin(\xii))^2$, whereas $\sin^{(k)}$ denotes the $k$-th derivative of $\sin$. We put $\sinc(\xii)=\xii^{-1}\sin(\xii)$.
\end{notation}
\begin{notation}
We use the notation $\<\xi\>=(1+\xii^2)^{\frac12}$, for $\xi\in \R^n$.
\end{notation}
\begin{notation}
For a given $b\in\R$, we say that $\phi\in\mathcal C^\infty(\R^n)$ is in $S^{-b}$ if
\[ |\partial_\xi^\alpha\phi(\xi)|\leq C_\alpha\,\<\xi\>^{-b-|\alpha|}. \]
To avoid confusion, we write $\S^{n-1}=\{\xi\in\R^n: \ \xii=1\}$ for the unit sphere in $\R^n$.
\end{notation}
\begin{notation}
Let $f\in\mathcal S'$ and $g(t,\cdot)\in\mathcal S'$ for any $t\geq0$. In this paper, $\mathscr{F}$ denotes the Fourier transform with respect to $x$ and we write $\hat f=\mathscr{F}f$ and $\hat g(t,\cdot)=\mathscr{F}g(t,\cdot)$. For $f\in L^1$,
\[ \hat f(\xi)=\int_{\R^n} e^{-ix\xi}\,f(x)\,dx. \]
Moreover, $\mathscr{F}^{-1}(f(\xi))$ and $\mathscr{F}^{-1}(g(t,\xi))$ denote the inverse Fourier transform with respect to $\xi$. When $\hat f\in L^1$,
\[ f(x)=(2\pi)^{-n}\,\int_{\R^n} e^{ix\xi}\,\hat f(\xi)\,d\xi. \]
With abuse of notation, we also write
\[ \|\mathscr{F}^{-1}(g(t,\xi))\|_{L^r} \]
to mean that the $L^r$ norm is taken with respect to the $x$ variable of the inverse Fourier transform of $\hat g(t,\cdot)$.
\end{notation}
\begin{notation}\label{not:chi}
In the following, we often use cut-off functions $\chi\in\mathcal C_c^\infty(\R^n)$, with $\chi=1$ in a neighborhood of the origin, and for any $t>0$ we put $\chi_t(\xi)=\chi (t^{\frac1\sigma}\xi)$. We make use of the localizing operators
\[ T_\chi f = \mathscr{F}^{-1}(\chi\,\hat f), \qquad T_{\chi_t} f = \mathscr{F}^{-1}(\chi_t\,\hat f). \]
\end{notation}
\begin{notation}[see~\cite{Ho}]\label{DefLpqspaces}
Let $1\leq p\leq q\leq\infty$. Then $M_p^q$ denotes the space of multipliers of type $(p,q)$ with norm defined by~\eqref{eq:Mpq}.
%
\end{notation}
\begin{notation}\label{not:H1}
In this paper, $H^1$ denotes the real $n$-dimensional Hardy space introduced by C. Fefferman and E. Stein~\cite{FS72}, equipped with norm~\eqref{eq:H1norm}; $M(H^1,L^q)$ and $M(H^1,H^1)$ denote the space of bounded multipliers from $H^1$ to $L^q$ or to itself, equipped with norms~\eqref{eq:H1multnorms}.
%
%
\end{notation}
\begin{notation}\label{not:Bessel spaces}
In this paper, $H^{\kappa,p}(\R^n)$, with $\kappa\geq0$ and $p\in[1,\infty]$, denotes the Bessel potential space
\[ H^{\kappa,p}(\R^n) = \{f\in L^p(\R^n): \quad \mathscr{F}^{-1}(\<\xi\>^\kappa\hat f)\in L^p(\R^n)\},  \]
equipped with norm
\[ \|f\|_{H^{\kappa,p}} = \|\mathscr{F}^{-1}(\<\xi\>^\kappa\hat f)\|_{L^p}. \]
Clearly, $H^{0,p}(\R^n)=L^p(\R^n)$. If $p\in(1,\infty)$ and $\kappa\in\N$, then $H^{\kappa,p}(\R^n)=W^{\kappa,p}(\R^n)$, the Sobolev space of $L^p(\R^n)$ functions with weak derivatives of order $\kappa$ in $L^p(\R^n)$.
\end{notation}


\section{Examples for dissipative wave equations}\label{sec:examples}

In this section, we provide examples of application of Theorems~\ref{thm:Wavedec} and~\ref{thm:Wavereg} for the wave equation, i.e. $L_{w^2}=-\Delta$ in~\eqref{eq:CPwave}, with different dissipative operators $A$, that is, we consider
\begin{equation}
\label{eq:CPwaveclassic}
\begin{cases}
u_{tt} -\Delta u + Au_t =0, \quad t>0,\, x\in \R^n, \\
u(0,x)=0, \\
u_t(0,x)=u_1(x).
\end{cases}
\end{equation}
Our results are indeed applicable to more general homogeneous hyperbolic equations and $\sigma$-evolution equations, but the peculiarities of the wave equation make the examples of particular interest, in our opinion. Therefore, for the sake of brevity, we focus on those examples.


The simplest model is when $A=(-\Delta)^{\frac\theta2}$, i.e., the case of wave equation with a so-called structural damping. In such a case, $a(\xi)=\xii^\theta$ is radial homogeneous. In recent years, low frequencies estimates in the case of so-called ``effective damping'' $\theta\in(0,1]$ have been well-investigated. 
In particular, the diffusion phenomenon has been proved in $L^p-L^q$ setting in~\cite{MN03, N03, HM} for $\theta=0$, and in~\cite{DAE2014} for $\theta\in(0,1)$ and $\sigma=1$. 
Similarly, high frequencies estimates in the case of the so-called ``noneffective damping'', $\theta\in[1,2]$, have been well-investigated. These estimates were relatively easy to be obtained due to the fact that oscillations are canceled in the overdamping regime, see~\eqref{eq:wavesolOD}. The study of $L^1-L^1$ estimates for these models goes back to~\cite{NR}.

In the case $\sigma\neq1$ and $\theta\in[\sigma,2\sigma]$, low frequencies estimates have been recently obtained by the authors~\cite{DAE21}, but the more difficult case of wave equation $\sigma=1$ remained open, so far.

\subsection{The wave equation with viscoelastic damping}\label{sec:visco}

When $A=-\Delta$, that is, $\theta=2$, the equation in~\eqref{eq:CPwaveclassic} reads as
\[ u_{tt} -\Delta u -\Delta u_t =0, \quad t>0,\, x\in \R^n, \]
it has been investigated in~\cite{Shi20} (see also~\cite{Ponce}) and it is related to the linearized Navier-Stokes equation in the following sense. Deriving with respect to $t$ the first equation in
\[ \begin{split}
& \rho_t+\nabla\cdot v = 0,\\
& v_t -\alpha\Delta v -\beta \nabla \nabla\cdot v +\nabla\rho=0,
\end{split} \]
and applying the divergence to the second one, it follows that the density $\rho$ satisfies the second order equation
\[ \rho_{tt} - (\alpha+\beta)\Delta\rho_t -\Delta\rho=0. \]
In~\cite[Theorem 2.1]{Shi20}, exception given for the trivial $L^p-L^q$ estimates when $p=2$ or $q=2$ (those estimates are trivial because the oscillations have no influence on the estimates), low frequencies $L^1-L^\infty$ estimates and $L^1-L^1$ estimates are obtained:
\begin{align}
\label{eq:Shiest1inf}
\|T_\chi u(t,\cdot)\|_{L^\infty}& \leq C\,(1+t)^{-\frac{3(n-1)}4}\,\|u_1\|_{L^1}, \\
\label{eq:Shiest1odd}
\|T_\chi u(t,\cdot)\|_{L^1}& \leq C\,(1+t)^{\frac{n+1}4}\,\|u_1\|_{L^1},\qquad \text{if $n\geq 3$ is odd,} \\
\label{eq:Shiest1even}
\|T_\chi u(t,\cdot)\|_{L^1}& \leq C\,(1+t)^{\frac{n+2}4}\,\|u_1\|_{L^1}, \qquad \text{if $n\geq2$ is even.}
\end{align}
Thanks to Theorem~\ref{thm:Wavedec} and Remark~\ref{Rem:noAw}, we find, for any $1\leq p\leq q\leq\infty$,
\[ \|T_\chi u(t,\cdot)\|_{L^q} \lesssim (1+t)^{1-n\left(\frac1p-\frac1q\right)+\frac{(d(p,q)-1)_+}2}\,\|u_1\|_{L^p}, \]
where $d$ is as in~\eqref{eq:dpqdef}, exception given for the special case $n=2$ and $(p,q)\in\{(1,2),(2,\infty)\}$, for which Proposition~\ref{prop:Wavedec2} gives us the estimate~\eqref{eq:Wavedec2}, where a logarithmic loss (sharp, see Lemma~\ref{lem:crucialwave2}) appears:
\begin{equation}\label{eq:L2crit}
\|T_\chi u(t,\cdot)\|_{L^q} \leq C (1+(\log t)_+)^{\frac12} \,\|u_1\|_{L^p}.
\end{equation}
In particular, when $p=q=1$ and $n\geq1$, due to $d=(n-1)/2$, we find
\[ \|T_\chi u(t,\cdot)\|_{L^1} \lesssim (1+t)^{\frac{n+1}4}\,\|u_1\|_{L^1}. \]
This latter is consistent with \eqref{eq:Shiest1odd}, while it improves the estimate in~\eqref{eq:Shiest1even}. For general $1\leq p\leq q\leq\infty$, our results are better than the obtained by interpolating the endpoints \eqref{eq:Shiest1inf} and~\eqref{eq:Shiest1odd} in~\cite{Shi20}. 

\subsection{The wave equation with noneffective dissipation}\label{sec:noneff}

As a variant of the viscoelastic damped wave equation, we may consider a wave equation with a noneffective damping, that is, $A=(-\Delta)^{\frac\theta2}$, with $\theta\in(1,2)$. In particular, if $1<p\leq q<\infty$, Theorem~\ref{thm:Wavedec} and Remark~\ref{Rem:noAw} provide the decay estimate
\begin{equation}\label{eq:neffest}
\|T_\chi u(t,\cdot)\|_{L^q} \leq C\,(1+t)^{1-n\left(\frac1p-\frac1q\right)+(d(p,q)-1)_+\left(1-\frac1\theta\right)}\,\|u_1\|_{L^p},
\end{equation}
exception given for the special case $n=2$ and $(p,q)\in\{(1,2),(2,\infty)\}$, for which Proposition~\ref{prop:Wavedec2} gives us estimate~\eqref{eq:L2crit}. These estimates are completely new. 

The wave equation with the so-called log-damping $a(\xi)=\log (1+\xii^\theta)$ has been recently considered in~\cite{CDI21,CI20}. Even if this damping is not homogeneous, due to $\log(1+\xii^\theta)\approx \xii^\theta$ as $\xi\to0$, the regime at low frequencies is of damped oscillations if $\theta>1$; hence, we obtain~\eqref{eq:neffest} applying Theorem~\ref{thm:Wavedec}. This result is also new for $q\neq2$.

\subsection{The wave equation with effective damping}\label{sec:eff}

When $A=(-\Delta)^{\frac\theta2}$, with $\theta\in(0,1)$, a ``double diffusion phenomenon'' holds, in the sense that the solution to~\eqref{eq:CPwaveclassic} may be written as the sum of two terms, whose asymptotic profiles as $t\to \infty$ are described by the solutions to the two diffusion problems \cite{DAE2014}. At low frequencies, due to the overdamping regime~\eqref{eq:wavesolOD}, decay estimates are easily obtained.

Theorem~\ref{thm:Wavereg} provides new, sharp, high frequencies $L^p-L^q$ estimates. In particular, they are consistent with the result obtained only in the case $p=q$ in~\cite{DAE2014nonlinear} when $(n-1)|1/p-1/2|\leq1$, that is, $d(p,p)\leq1$ (even restricting only to $L^p-L^p$ estimates, Theorem~\ref{thm:Wavereg} provides estimate also when $d(p,p)>1$). In particular, taking $s=0$ in Theorem~\ref{thm:Wavereg} and taking into account of Remark~\ref{Rem:noAw}, we get:
\[ \|(\id-T_\chi)u(t,\cdot)\|_{L^q} \leq C\,t^{1-n\left(\frac1p-\frac1q\right)-(d(p,q)-1)_+\left(\frac1{\theta}-1\right)}\,e^{-ct}\,\|u_1\|_{L^p},\]
for any $t>0$, exception given for the case $n=2$ and $(p,q)\in\{(1,2),(2,\infty)\}$, for which Proposition~\ref{prop:Wavereg2} gives us estimate~\eqref{eq:Wavereg2}, where a logarithmic singularity (sharp, see Lemma~\ref{lem:crucialwave2}) appears:
\[ \|(\id-T_\chi)\,u(t,\cdot)\|_{L^2} \leq C\,(1+(-\log t)_+)^{\frac12}\,e^{-ct}\,\|u_1\|_{L^1},\quad t>0.\]

\subsection{The double dispersion equation}\label{sec:dde}

The equation
\[ u_{tt}-\Delta u_{tt} -\Delta u + \Delta^2 u + Bu_t=0, \]
is called (damped) double dispersion equation and has been studied by several authors, since it is related to a generalized Boussinesq equation (see~\cite{Chen, Polat, Wang} and the references therein). This model is also called damped plate equation with rotational inertia~\cite{CLI,LC,S-K}.

By applying the Bessel potential $(\id-\Delta)^{-1}$, the equation reduces to the one in~\eqref{eq:CPwaveclassic}, with~$A=(\id-\Delta)^{-1}B$. If $B=(-\Delta)^{\frac\theta2}$, for some $\theta$, we find that
\[ a(\xi)= \frac{\xii^\theta}{1+\xii^2}, \]
in particular, is not homogeneous. This model is also of interest because we have a damped oscillations regime at the both low and high frequencies. Letting $\theta_0=\theta$ and $\theta_1=\theta-2$, we find that Theorem~\ref{thm:Wavedec} is applicable for $\theta>1$ and Theorem~\ref{thm:Wavereg} is applicable for $\theta\in(2,3)$. The case $\theta=2$ is of particular interest; Theorem~\ref{thm:Wavedec} and Remark~\ref{Rem:noAw} provide the following estimates, that are new as far as we know:
\[
\|T_\chi\,u(t,\cdot)\|_{L^q} \leq C\,(1+t)^{1-\frac{n}\sigma\left(\frac1p-\frac1q\right)+\frac{(d-1)_+}2}\,\|u_1\|_{L^p}, \quad t\geq0,
\]
exception given, as usual, for the case $n=2$ and $(p,q)\in\{(1,2),(2,\infty)\}$, for which Proposition~\ref{prop:Wavedec2} gives us~\eqref{eq:L2crit}.

\section{Multiplier theorems}\label{sec:multipliers}

We recall some multipliers properties that we will use in this paper.

We stress that $M_p^q=\{0\}$ when $1\leq q<p<\infty$ and that $M_p^q$ contains multipliers that are distributions with positive order if, and only if, $1\leq p<2<q\leq\infty$, see~\cite[Theorem 1.6]{Ho}; indeed, $M_p^q\subset L_\loc^{q'}$ if $q\leq2$, and $M_1^q=L^q$ for any $q\in(1,\infty]$, while $M_1^1$ is the set of bounded measures, see~\cite[Theorem 1.4]{Ho}. %
In this paper, multipliers will be functions, though, since they are related to the Fourier transform of the fundamental solutions to~\eqref{eq:CPSch} and~\eqref{eq:CPwave}.

It holds
\[ M_{p_1}^{p_1} \subset M_{p_2}^{p_2}\subset M_2^2=L^\infty, \qquad 1\leq p_1\leq p_2\leq 2, \]
and (see~\cite[Theorem 1.4]{Ho})
\begin{equation}\label{eq:M1p}\begin{split}
M_1^p
    & =\{m\in\mathcal S': \ \mathscr{F}^{-1}(m)\in L^{p}\},\ p\in(1,\infty]; \\
M_1^1
    & =\{m\in\mathcal S': \ \text{$\mathscr{F}^{-1}(m)$ is a bounded measures}\},
\end{split}\end{equation}
in particular, $M_1^2=L^2$, by Plancherel's theorem. By Young inequality, it holds
\[ \{m\in \mathcal S': \ \mathscr{F}^{-1}(m)\in L^r\} \subset M_p^q, \quad 1-\frac1r = \frac1p-\frac1q. \]
Let $H^1$ denote the real Hardy space (see Notation~\ref{not:H1}).Then, it holds (see~\cite[Theorem 3.3]{Miy81}):
\begin{equation}\label{eq:M1H1}
\|m\|_{M(L^1,H^1)}=\|\mathscr{F}^{-1}m\|_{H^1}.
\end{equation}
On the other hand (see~\cite{FS72}, see also~\cite[Theorem 3.4]{Miy81}):
\[ \|m\|_{M(H^1,L^1)}=\|m\|_{M(H^1,H^1)}.\]
We stress that
\begin{equation}\label{eq:H1L1order}
M(L^1,H^1)\subset M_1^1\subset M(H^1,H^1)=M(H^1,L^1).
\end{equation}
By Hardy-Littlewood-Sobolev inequality, if $1<p<\infty$ and $0<a<n/p$,
\begin{equation}\label{eq:HLS}
\xii^{-a} \in M_p^{p^*}, \quad \frac1{p^*}=\frac1p-\frac{a}n.
\end{equation}
The above property holds for $p=1$ if $M_1^{1^*}$ is replaced by $M(H^1,L^{1^*})$ (see~\cite[Theorem F]{Miy81}).

By Mikhlin-H\"ormander theorem (see also~\cite{Grafakos2021} for several variants of this result), for any $p\in(1,\infty)$, $M_p^p$ contains the space of $\mathcal C^{1+[n/2]}(\R^n\setminus\{0\})$ multipliers (here $[\cdot]$ denotes the floor function) such that
\[ |\partial_\xi^\gamma m(\xi)|\leq C\,\xii^{-|\gamma|}, \quad |\gamma|\leq 1+[n/2]. \]
This space is also contained in $M(H^1,H^1)$ (see~\cite[Theorem E]{Miy81}). The Mikhlin-H\"ormander theorem may be used to easily prove that the multipliers studied in Lemmas~\ref{lem:exp}, \ref{lem:exp2} and \ref{lem:exp3} are in $M_p^p$ when $p\in(1,\infty)$, and that the regularity of $a(\xi)$ may be reduced to $\mathcal C^{1+[n/2]}(\R^n\setminus\{0\})$, asking that~\eqref{eq:oscillationslow} and~\eqref{eq:oscillationshigh} hold only for  $|\gamma|\leq 1+[n/2]$. However, this approach cannot be applied in the endpoints $p=1,\infty$.

To prove Theorem~\ref{thm:w}, we will also make use of the following theorem for $M(H^1,H^1)$.
\begin{Thm}\label{thm:MHardy}[Theorem G in~\cite{Miy81}, see also Theorems 1 and 2 in~\cite{Miy80}]
Assume that $m\in\mathcal C^{1+[n/2]}(\R^n\setminus\{0\})$ are such that
\[ |\partial_\xi^\alpha m(t,\xi)| \leq \<\xi\>^{-\frac{n}2\sigma} \,\xii^{-|\alpha|} ((1+t)\,\<\xi\>^{\sigma-1})^{|\alpha|},\quad |\alpha|\leq 1+[n/2], \]
for $t\geq0$ and $\sigma>0$. Then $m(t,\xi)\in M(H^1,H^1)$ and
\[ \|m(t,\xi)\|_{M(H^1,H^1)} \leq (1+t)^{\frac{n}2}. \]
\end{Thm}
By the homogeneous property of the Fourier transform,
\[ \|m(t\xi)\|_{M_p^q} = t^{-n\left(\frac1p-\frac1q\right)}\,\|m(\xi)\|_{M_p^q}. \]
This property clearly also holds when $p=1$ if $M_1^q$ is replaced by $M(H^1,L^q)$, or by $M(H^1,H^1)$ or $M(L^1,H^1)$ when $q=1$ (see~\cite[Lemma 3.3]{Miy81}).

To deal with regularity of initial data, we are also going to use that
\begin{equation}\label{eq:M11reg}
\xii^\beta\,\<\xi\>^{-\beta}\in M_1^1,\quad \forall \beta\geq0,
\end{equation}
since the multiplier is the Fourier transform of a bounded measure (see, for instance, \cite[\textsection3.2, Lemma 2]{Stein70}).

\bigskip

We also recall the following result (see \cite[Theorem 7]{FS72}), which may be applied to obtain estimates for (fractional) derivatives in $x$ of the solution to~\eqref{eq:CPSch} and~\eqref{eq:CPwave}. If $m\in M(H^1,H^1)$ and $|m(\xi)|\leq C\xii^{-\delta}$ for some $\delta>0$, then $\xii^\theta m \in M_p^p$ for any $\theta\in[0,\delta]$ and $p\in(1,\infty)$ such that
\[ \left|\frac1p-\frac12\right|\leq \frac12-\frac\theta{2\delta}. \]
In general, our estimates in this paper may be easily modified to take into account of the presence of (fractional) derivatives of $u$, or of time derivatives of $u$, but we avoid this study for the sake of brevity.

\bigskip

We are now interested in estimates for diffusive multipliers in the form $e^{-ta(\xi)}$ and their perturbations $b(\xi)\,e^{-ta(\xi)}$.

In order to deal with multipliers in $M_1^1$, we cannot rely on Mikhlin-H\"ormander multiplier theorem. A widely used alternative is the integration by parts method. We use the integration by parts method to prove several results for parameter-dependent multipliers of the form $e^{-ta(\xi)}$, whose proof we postpone to \textsection\ref{sec:multproofs} for the ease of reading.

The first basic result is the following.
\begin{lemma}\label{lem:exp}
Let $a\in\mathcal C^{n+1}(\R^n\setminus\{0\})$ verifies
\[ a(\xi)\geq c\,\xii^{\theta},\quad |\partial_\xi^\gamma a(\xi)|\leq C\,\xii^{\theta-|\gamma|}, \qquad \text{for $1\leq|\gamma|\leq n+1$,}\]
for some~$c>0$, $C>0$, $\theta>0$. Then $\mathscr{F}^{-1} (e^{-ta(\xi)})$ is in $L^r$ and
\[ \|\mathscr{F}^{-1} (e^{-ta(\xi)})\|_{L^r} \lesssim t^{-\frac{n}\theta\left(1-\frac1r\right)}, \]
for any $r\in[1,\infty]$.
\end{lemma}
We stress that if $a(\xi)$ is radial homogeneous, that is, $a(\xi)=c\xii^\theta$, Lemma~\ref{lem:exp} follows as a consequence of the homogeneity property of $\mathscr{F}$ and by the asymptotic behavior~\cite{BG}:
\begin{equation}\label{eq:asymp}
\lim_{|x|\to\infty} |x|^{n+\theta}\,\mathscr{F}^{-1} (e^{-\xii^\theta})=C_{n,\theta}.
\end{equation}
We also have the following generalizations of Lemma~\ref{lem:exp}.
\begin{lemma}\label{lem:exp2}
Let $a\in\mathcal C^{n+1}(\R^n\setminus\{0\})$ verifies
\[ a(\xi)\geq c\,\xii^{\theta},\quad |\partial_\xi^\gamma a(\xi)|\leq C\,\xii^{\theta-|\gamma|},\qquad \text{for $1\leq|\gamma|\leq n+1$,} \]
for some~$c>0$, $C>0$, $\theta>0$, and let $b\in\mathcal C^{n+1}(\R^n\setminus\{0\})$ verify
\[ |\partial_\xi^\gamma b(\xi)|\leq C\,\xii^{\eta-|\gamma|},\quad \qquad \text{for $0\leq|\gamma|\leq n+1$,}\]
for some $\eta>-n$. Then
\[ \|\mathscr{F}^{-1} (b(\xi)\,e^{-ta(\xi)})\|_{L^r}\lesssim t^{-\frac{n}\theta\left(1-\frac1r\right)-\frac\eta\theta}, \]
for any $r\in[1,\infty]$, if~$\eta>0$ and for any $r\in(1,\infty]$ such that
\begin{equation}\label{eq:reta}
n\left(1-\frac1r\right)>-\eta,
\end{equation}
if $\eta\in(-n,0]$.
\end{lemma}
The following consequence of Lemma~\ref{lem:exp2} is straightforward, but not very standard, since it involves the real Hardy space $H^1$. Its use is crucial to obtain sharp estimates with $p=1$ or $q=\infty$ in Theorems~\ref{thm:Schdec}, \ref{thm:Schreg}, \ref{thm:Wavedec} and~\ref{thm:Wavereg}.
\begin{Cor}\label{cor:H1}
Assume that $\eta>0$ in Lemma~\ref{lem:exp2}. Then the $L^1$ estimate may be improved to
\[ \|\mathscr{F}^{-1} (b(\xi)\,e^{-ta(\xi)})\|_{H^1} \lesssim t^{-\frac\eta\theta},\]
where $H^1$ denotes the real Hardy space.
\end{Cor}
\begin{proof}
In view of~\eqref{eq:H1norm}, it is sufficient to apply Lemma~\ref{lem:exp2} with $r=1$ and $b_j=b(\xi)\xi_j/\xii$ in place of $b$, since
\[ \|\mathscr{F}^{-1} (b(\xi)\,e^{-ta(\xi)})\|_{H^1} = \|\mathscr{F}^{-1} (b(\xi)\,e^{-ta(\xi)})\|_{L^1} + \sum_{j=0}^n\|\mathscr{F}^{-1} (b_j(\xi)\,e^{-ta(\xi)})\|_{L^1}.\]
\end{proof}
The next generalization of Lemma~\ref{lem:exp} is of interest to treat several critical cases and we will only use it at low frequencies.
\begin{lemma}\label{lem:exp3}
Let $a\in\mathcal C^{n+1}(\R^n\setminus\{0\})$ verifies
\[ a(\xi)\geq c\,\xii^{\theta},\quad |\partial_\xi^\gamma a(\xi)|\leq C\,\xii^{\theta-|\gamma|},\qquad \text{for $1\leq|\gamma|\leq n+1$,} \]
for some~$c>0$, $C>0$, $\theta>0$, and let $b\in\mathcal C^{n+1}(\R^n\setminus\{0\})$ verifies
\[ |b(\xi)|\leq C,\quad \text{and}\qquad |\partial_\xi^\gamma b(\xi)|\leq C\,\xii^{\eps-|\gamma|}, \qquad \text{for $1\leq|\gamma|\leq n+1$,} \]
for some $\eps>0$. Then
\[ \|\mathscr{F}^{-1} (b(\xi)\,e^{-ta(\xi)})\|_{L^r} \lesssim t^{-\frac{n}\theta\left(1-\frac1r\right)},\quad \text{for $t\geq1$,} \]
for any $r\in[1,\infty]$.
\end{lemma}
%
Finally, we need the following straightforward result.
\begin{lemma}\label{lem:hexp3}
Let $\phi\in S^{-\eps}$ for some $\eps>0$ or, more in general, assume that $\phi\in\mathcal C^{n+1}(\R^n)$ verifies
\begin{equation}\label{eq:strongder}
|\partial_\xi^\gamma \phi(\xi)|\leq C\,\<\xi\>^{-\eps-|\gamma|},\qquad \text{for $|\gamma|=n-1,n,n+1$,}
\end{equation}
for some $\eps\in(0,1)$. Then $\mathscr{F}^{-1} \phi\in L^1$. Moreover, $\phi\in M(L^1,H^1)$ if $\phi$ vanishes in a neighborhood of the origin.
\end{lemma}
%


\section{Proof of Theorems~\ref{thm:w}--\ref{thm:A0w}}\label{sec:w}

We preliminarily notice that, thanks to the homogeneity assumption of $w$, letting $\xi'=\xi/\xii$, we have the following:
\begin{align}
\label{eq:MHhom}
|\partial_\xi^\gamma w(\xi)| & = \xii^{\sigma-|\gamma|} |\partial_\xi^\gamma w(\xi')| \leq C_\gamma\,\xii^{\sigma-|\gamma|},\\
\label{eq:poshom}
w(\xi) & =\xii^\sigma\,w(\xi')\geq c\xii^\sigma, \qquad c=\min_{|\xi|=1} w(\xi)>0.
\end{align}

We now prove Theorem~\ref{thm:w}.
\begin{proof}[Proof of Theorem~\ref{thm:w}]
Since Theorem~\ref{thm:w} is already proved for $\sigma=1$, we assume $\sigma\neq1$.

We prove~\eqref{eq:wbasic} using the complex interpolation (for instance, see \cite{Stein71} Chapter 5 \textsection 4; see also \cite[\textsection1.3.3]{Grafakos14book}). We first consider the case $q=p\in(1,\infty)$ in~\eqref{eq:wbasic}, and we assume without restriction that $b<n/2$. For $s\in [0,1]+i\R$, we define the multiplier
\[ m_s= e^{(s-\theta)^2}\,\phi\,\<\xi\>^{\gamma(s)\sigma}\,e^{\pm iw(\xi)}, \]
%
%
%
%
where
\[ \gamma(s)=b-s\frac{n}2,\quad \theta=\frac{2b}n. \]
We notice that $\gamma(\theta)=0$, so that $m_\theta=\phi\,e^{\pm iw(\xi)}$. 

When $\Re s=0$ and $\Im s=t\in\R$,
\[|m_{it}(\xi)|=|e^{(it-\theta)^2} \phi(\xi)\,\<\xi\>^{\gamma(it)\sigma}|\leq C\,|e^{(it-\theta)^2}|\leq e^{\theta^2}=C_1,\]
with $C_1$ independent on $\Im s$, so that
\[ \|m_{it}\|_{M_2^2}=\|m_{it}\|_{L^\infty}\leq C_1.\]
When $\Re s=1$ and $\Im s=t\in\R$, $\phi(\xi)\,\<\xi\>^{\gamma(1+it)\sigma} \in S^{-\frac{n}2\sigma}$; applying Theorem~\ref{thm:MHardy}, due to
\[ |\partial_\xi^\alpha m_{1+it}(\xi)|\leq C\,|e^{(1+it-\theta)^2}|\,(1+|t|)^{|\alpha|}\, |\xi|^{-|\alpha|}\,\<\xi\>^{|\alpha|\sigma-\frac{n}2\sigma},\quad |\alpha|\leq \frac{n}2+1, \]
where we used~\eqref{eq:MHhom}, we find that
\begin{equation}\label{eq:H1H1int}
\|m_{1+it}\|_{M(H^1,H^1)}\leq C_1\,e^{(1-\theta)^2-t^2}\,(1+|t|)^{\frac{n}2}\leq C_2.
\end{equation}
We now fix $p\in(1,2]$ with
\[ \frac1p=\frac12 + \frac{b}n\]
and we use complex interpolation (see \cite[\textsection5, Corollary 1]{FS72}), so that $m_\theta\in M_p^p$. Thanks to the  duality argument we conclude the first part of~\eqref{eq:wbasic} for $q=p\in(1,\infty)$.

If $b\geq n/2$, it is sufficient to use $H^1-H^1$ estimate~\eqref{eq:H1H1int}, with $t=0$, to derive the second part of~\eqref{eq:wbasic}, with $q=1$ (we recall that $M(H^1,H^1)=M(H^1,L^1)$, see~\eqref{eq:H1L1order}).

We now prove~\eqref{eq:endpoint}. Thanks to the Hessian assumption, we may assume without restriction that for any given $x$, there is at most one critical point for $x\xi\pm w(\xi)$ in the support of $\phi$. Let $\xi_0$ be this critical point, that is, $x=\mp \nabla w(\xi_0)$, and $\xi_0=\xi_0'\rho$, with $\rho=|\xi_0|$. Fix $\supp\phi\subset B(\xi_0,\delta)$, for some $\delta>0$ and $\delta=\rho\delta'$. By the change of variable $\xi\mapsto\rho\xi$, we obtain
\[ \begin{split}
K(x)
    & = \mathscr{F}^{-1}(e^{\pm iw(\xi)}\phi) =(2\pi)^{-n}\,\int_{B(\xi_0,\delta)} e^{\pm i(w(\xi)-\xi\nabla w(\xi_0))}\,\phi(\xi)\,d\xi\\
    & =(2\pi)^{-n}\,\rho^n\,\int_{B(\xi_0',\delta')} e^{\pm i(w(\rho\xi)-\rho\xi\nabla w(\xi_0))}\,\phi(\rho\xi)\,d\xi\\
    & =(2\pi)^{-n}\,\rho^n\,\int_{B(\xi_0',\delta')} e^{\pm i\rho^\sigma (w(\xi)-\xi\nabla w(\xi_0'))}\,\phi(\rho\xi)\,d\xi,
\end{split} \]
where in the last equality we used the homogeneity of $w$ and $\nabla w$. Due to the fact that $H_w$ is nonsingular, we may apply Littman's lemma (see, for instance, Proposition 6 in VIII in~\cite{Stein93}) and conclude that
\[ |K(x)| \sim (2\pi)^{-n}\rho^n\,(1+\rho^\sigma)^{-\frac{n}2}\,|\phi(\rho\xi_0')| \leq C, \]
thanks to assumption
\[ b\geq d(1,\infty)=n\frac{2-\sigma}{2\sigma}\,.\]
In absence of critical points, $|K(x)|\leq C$ also holds by Littman's lemma; hence, we obtained~\eqref{eq:endpoint}. By interpolation, proceeding as in the first step, this also concludes the proof of~\eqref{eq:wbasic}, with $p<q$.

To prove~\eqref{eq:w1}, we write $\phi=\phi_d\,\phi_\varepsilon$, where $\phi_d\in S^{-\sigma d}$ and $\phi_\varepsilon\in S^{-\sigma\varepsilon}$, with $\varepsilon=b-d>0$. The proof follows since $\phi_\varepsilon\in M(L^1,H^1)$ by Lemma~\ref{lem:hexp3}.

Without the maximal rank assumption, the proof follows for $1\leq p\leq q\leq2$, interpolating $L^2-L^2$, $H^1-H^1$ and $H^1-L^2$ estimates; the latter can be easily derived by using the Hardy-Littlewood-Sobolev inequality in real Hardy spaces~\eqref{eq:HLS}.

This concludes the proof of Theorem~\ref{thm:w}.
\end{proof}
%
%
Having in mind to use the homogeneity of $w(\xi)$ in $e^{i\pm w(\xi)t}$, for a given $\chi\in\mathcal C_c^\infty$ with $\chi=1$ in a neighborhood of the origin, and for any $t>0$ we use the cut-off function $\chi_t(\xi)=\chi(t^{\frac1\sigma}\xi)$ (see Notation~\ref{not:chi}). It holds:
\begin{equation}\label{eq:osclow}
\|\chi_t e^{\pm iw(\xi)t}\|_{M_p^q}
    = t^{-\frac{n}\sigma\left(\frac1p-\frac1q\right)}\,\|\chi e^{\pm iw(\xi)}\|_{M_p^q}= C\,t^{-\frac{n}\sigma\left(\frac1p-\frac1q\right)}
\end{equation}
for any $1\leq p\leq q\leq\infty$. Indeed, the fact that $\chi e^{\pm iw(\xi)}\in M_p^q$ for any $1\leq p\leq q\leq\infty$, in \eqref{eq:osclow}, is pretty standard. Thanks to~\eqref{eq:MHhom}, for any $\xi\neq0$ it holds
\[ |\partial_\xi^\alpha (\chi e^{\pm iw(\xi)})| \leq C\,\xii^{\sigma-|\alpha|}, \qquad |\alpha|=1,\ldots,n+1. \]
For instance, it is sufficient to apply Lemma 8.5 in~\cite{DAE21} to obtain that $\chi e^{\pm iw(\xi)}\in M_p^q$. In the case $1<p\leq q<\infty$, it is sufficient to combine Mikhlin-H\"ormander theorem and Hardy-Littlewood inequality.

On the other hand, $\phi = (1-\chi)\,\xii^{-\beta\sigma}\in S^{-\beta\sigma}$, vanishes in a neighborhood of the origin, so that we may apply Theorem~\ref{thm:w}. In particular,
\begin{align}
\label{eq:oschi}
\begin{split}
\|(1-\chi_t)\,\xii^{-\beta\sigma}\, e^{\pm iw(\xi)t}\|_{M_p^q}
    & = t^{\beta-\frac{n}\sigma\left(\frac1p-\frac1q\right)}\,\|(1-\chi)\,\xii^{-\beta\sigma}\, e^{\pm iw(\xi)}\|_{M_p^q} \\
    & = C\,t^{\beta-\frac{n}\sigma\left(\frac1p-\frac1q\right)},
\end{split}
\intertext{for~$\beta\geq d(p,q)$ if $1<p\leq q<\infty$ (and $(p,q)=(1,\infty)$ if $\sigma\neq1$), or for~$\beta>d(p,q)$, otherwise. Moreover,}
\label{eq:oschiH1}
\begin{split}
\|(1-\chi_t)\,\xii^{-\beta\sigma}\, e^{\pm iw(\xi)t}\|_{M(H^1,L^q)}
    & = t^{\beta-\frac{n}\sigma\left(\frac1p-\frac1q\right)}\,\|(1-\chi)\,\xii^{-\beta\sigma}\, e^{\pm iw(\xi)}\|_{M(H^1,L^q)} \\
    & = C\,t^{\beta-\frac{n}\sigma\left(\frac1p-\frac1q\right)},\quad \text{for~$\beta\geq d(1,q)$.}
\end{split}
\end{align}
%
We may now prove Theorem~\ref{thm:A0}.
\begin{proof}[Proof of Theorem~\ref{thm:A0}]
%
%
Let $\hat u(t,\xi) = e^{\pm itw(\xi)}\,\hat u_0(\xi)$. On the one hand, we may estimate
\[ \|T_{\chi_t}u(t,\cdot)\|_{L^q}\leq \big\|\chi_t\,e^{\pm iw(\xi)t}\|_{M_p^q}\,\|u_0\|_{L^p}\leq C\,t^{-\frac{n}\sigma\left(\frac1p-\frac1q\right)}\|u_0\|_{L^p},\]
where we used~\eqref{eq:osclow}; on the other hand, recalling that $w(\xi)t=w(t^{\frac1\sigma}\xi)$, thanks to~\eqref{eq:M11reg} and~\eqref{eq:oschi}, we may estimate
\[ \begin{split}
\|(\id-T_{\chi_t})u(t,\cdot)\|_{L^q}
    & \leq \big\|(1-\chi_t)\,\xii^{-s\sigma}\,e^{\pm iw(t^{\frac1\sigma}\xi)}\|_{M_p^q}\,\|\xii^{s\sigma}\,\<\xi\>^{-s\sigma}\|_{M_q^q}\,\|u_0\|_{H^{s\sigma,p}}\\
    & = t^{s-\frac{n}\sigma\left(\frac1p-\frac1q\right)}\,\big\|(1-\chi_1)\,\xii^{-s\sigma}\,e^{\pm iw(\xi)}\|_{M_p^q}\,\|\xii^{s\sigma}\,\<\xi\>^{-s\sigma}\|_{M_q^q}\,\|u_0\|_{H^{s\sigma,p}}\\
    & \leq C\,t^{s-\frac{n}\sigma\left(\frac1p-\frac1q\right)}\|u_0\|_{H^{s\sigma,p}}.
\end{split}\]
This concludes the proof.
\end{proof}
The proof of Theorem~\ref{thm:A0w} is relatively similar to the proof of Theorem~\ref{thm:A0}. Similarly to~\eqref{eq:osclow}, it holds:
\begin{equation}
\label{eq:waveosclow}
\|\chi_t \sinc (w(\xi)t)\|_{M_p^q}= t^{-\frac{n}\sigma\left(\frac1p-\frac1q\right)}\,\|\chi \sinc(w(\xi))\|_{M_p^q}= C\,t^{-\frac{n}\sigma\left(\frac1p-\frac1q\right)},
\end{equation}
for any $1\leq p\leq q\leq\infty$. Indeed,
\[ |\partial_\xi^\alpha \sinc(w(\xi))| \leq C\,\xii^{2\sigma-|\alpha|}, \qquad |\alpha|=1,\ldots,n+1, \]
and it is sufficient to apply Lemma 8.5 in~\cite{DAE21} to obtain that $\chi \sinc(w(\xi))\in M_p^q$. Moreover, we are going to use~\eqref{eq:oschi} and~\eqref{eq:oschiH1}.
\begin{proof}[Proof of Theorem~\ref{thm:A0w}]
Let $u$ be as in~\eqref{eq:uwave}. On the one hand, we may estimate
\[ \|T_{\chi_t}u(t,\cdot)\|_{L^q}\leq t\,\big\|\chi_t\,\sinc(w(\xi)t)\|_{M_p^q}\,\|u_1\|_{L^p}\leq C\,t^{1-\frac{n}\sigma\left(\frac1p-\frac1q\right)}\|u_1\|_{L^p},\]
where we used~\eqref{eq:waveosclow}; on the other hand, recalling that $w(\xi)t=w(t^{\frac1\sigma}\xi)$, thanks to~\eqref{eq:M11reg} and~\eqref{eq:oschi}, we may estimate
\[ \begin{split}
\|(\id-T_{\chi_t})u(t,\cdot)\|_{L^q}
    & \leq \Big\|(1-\chi_t)\,\xii^{-s\sigma}\,\frac{e^{iw(t^{\frac1\sigma}\xi)}-e^{-iw(t^{\frac1\sigma}\xi)}}{2w(\xi)}\Big\|_{M_p^q}\,\|\xii^{s\sigma}\,\<\xi\>^{-s\sigma}\|_{M_q^q}\,\|u_1\|_{H^{s\sigma,p}}\\
    & \leq C\,t^{s+1-\frac{n}\sigma\left(\frac1p-\frac1q\right)}\|u_1\|_{H^{s\sigma,p}}.
\end{split}\]
When $w=c\xii$, in the two special cases $n=1$ and $(p,q)=(1,\infty)$, $n=3$ and $p=q=1$, there is no need to split in low and high frequencies, since
\[ \begin{split}
\|u(t,\cdot)\|_{L^\infty}
    & \leq t\,\|\sinc (c\xii t)\|_{M_1^\infty}\,\|u_1\|_{L^1} \\
    & = c^{-1}\,\|\sinc \xii \|_{M_1^\infty}\,\|u_1\|_{L^1}\leq C\,\|u_1\|_{L^1},\qquad n=1,\\
\|u(t,\cdot)\|_{L^p}
    & \leq t\,\|\sinc (c\xii t)\|_{M_1^1}\,\|u_1\|_{L^p} \\
    & = t\,\|\sinc \xii \|_{M_1^1}\,\|u_1\|_{L^p}=Ct\,\|u_1\|_{L^p},\qquad n=3.
\end{split}\]
As well known, this unique peculiarity is due to the fact that in space dimension $n=1$,
\[ \sinc\xii = \sinc \xi = \frac12\mathscr{F}(\chi_{[-1,1]}), \]
where $\chi_{[-1,1]}$ is the indicator function of the interval $[-1,1]$, in particular, it is in $L^\infty$ (this is also known as d'Alembert formula), while in space dimension $n=3$,
\[ \sinc\xii = \mathscr{F}(\delta_1), \]
where $\delta_1$ is the measure that gives the average of a function on the unit sphere; in particular, $\sinc\xii\in M_1^1=\{\text{bounded measures}\}$.

This concludes the proof.
\end{proof}

\section{Proof of Theorems~\ref{thm:Schdec} and~\ref{thm:Schreg}}\label{sec:Sch}

In this Section, we prove Theorems~\ref{thm:Schdec} and~\ref{thm:Schreg}, using Theorem~\ref{thm:w} and the results for the perturbation multiplier collected in Lemmas~\ref{lem:exp}, \ref{lem:exp2}, \ref{lem:exp3} and \ref{lem:hexp3}.

We first prove Theorem~\ref{thm:Schdec}.
\begin{proof}[Proof of Theorem~\ref{thm:Schdec}]
We want to estimate the following term:
\[ \|T_\chi u(t,\cdot)\|_{L^q}\leq \big\|e^{\pm iw(\xi)t}\,\chi \,e^{-a(\xi)t}\|_{M_p^q}\,\|u_0\|_{L^p}. \]
Since $\chi$ is compactly supported, the uniform estimate
\[ \big\|e^{\pm iw(\xi)t}\,\chi \,e^{-a(\xi)t}\|_{M_p^q}\leq C\,,\qquad t\in[0,1],\]
is trivial. Indeed, due to
\[ |\partial_\xi^\gamma (e^{\pm iw(\xi)t} \,e^{-a(\xi)t})| \leq C\,\xii^{\min\{\sigma,\theta_0\}-|\gamma|}\,, |\gamma|=1,\ldots,n+1, \]
for $\xi\in\supp\chi$ and for $t\in[0,1]$, it is sufficient to apply Lemma 8.5 in~\cite{DAE21}. %
%

Now let $t\geq1$. On the one hand, we get
\[ \begin{split}
\big\|\chi_t e^{\pm iw(\xi)t}\,\chi\,e^{-a(\xi)t}\|_{M_p^q}
    & \leq \|\chi_t e^{\pm iw(\xi)t}\|_{M_p^q}\,\|\chi\,e^{-a(\xi)t}\|_{M_q^q} \\
    & \leq C\,t^{-\frac{n}\sigma\left(\frac1p-\frac1q\right)},
\end{split} \]
where we used~\eqref{eq:osclow} and Lemma~\ref{lem:exp3} with $b=\chi$ (since $\chi=1$ in a neighborhood of the origin, we may take $\varepsilon>0$ as we want); on the other hand, if $1<p\leq q<\infty$, we estimate
\[ \begin{split}
\big\|(1-\chi_t) e^{\pm iw(\xi)t}\,\chi\,e^{-a(\xi)t}\big\|_{M_p^q}
    & \leq \|(1-\chi_t)\xii^{-\beta\sigma} e^{\pm iw(\xi)t}\|_{M_p^q}\,\|\chi\,\xii^{\beta\sigma}\,e^{-a(\xi)t}\|_{M_q^q} \\
    & \leq C\,t^{\beta-\frac{n}\sigma\left(\frac1p-\frac1q\right)-\beta\,\frac{\sigma}{\theta_0}},
\end{split} \]
where we used~\eqref{eq:oschi} with $\beta=d(p,q)$ and Lemma~\ref{lem:exp2} with $b=\chi\,\xii^{\beta\sigma}$. If $p=1$, $q\in[1,\infty]$, we use the real Hardy space $H^1$; thanks to~\eqref{eq:oschiH1}, and to Corollary~\ref{cor:H1}, we get the following estimate
\[ \begin{split}
& \big\|(1-\chi_t) e^{\pm iw(\xi)t}\,\chi\,e^{-a(\xi)t}\big\|_{M_1^q}\\
    & \qquad \leq \|(1-\chi_t)\xii^{-\beta\sigma} e^{\pm iw(\xi)t}\|_{M(H^1,L^q)}\,\|\chi\,\xii^{\beta\sigma}\,e^{-a(\xi)t}\|_{M(L^1,H^1)} \\
    & \qquad \leq C\,t^{\beta-\frac{n}\sigma\left(1-\frac1q\right)-\beta\,\frac{\sigma}{\theta_0}},
\end{split} \]
where $\beta=d(1,q)$. By duality, the same estimate holds for $q=\infty$ and $p\in[1,\infty]$. We notice that
\[ \beta-\beta\,\frac{\sigma}{\theta_0} \geq 0 \iff \theta_0\geq \sigma. \]
Therefore, %
%
%
%
%
we obtain the desired estimates and this concludes the proof.
\end{proof}
We now prove Theorem~\ref{thm:Schreg}.
\begin{proof}[Proof of Theorem~\ref{thm:Schreg}]
We want to estimate the following term:
\[ \|(\id-T_\chi) u(t,\cdot)\|_{L^q}=\big\|\mathscr{F}^{-1}\big(e^{\pm iw(\xi)t}\,(1-\chi)\,e^{-a(\xi)t}\,\hat u_0\big)\big\|_{L^q}. \]
First, let $\theta_1\in(0,\sigma)$. On the one hand,
\[ \begin{split}
& \big\|\mathscr{F}^{-1}\big(\chi_t e^{\pm iw(\xi)t}\,(1-\chi)\,e^{-a(\xi)t}\,\hat u_0\big)\big\|_{L^q}\\
& \qquad \leq \|\chi_t e^{\pm iw(\xi)t}\|_{M_p^q}\,\|(1-\chi)\,e^{-a(\xi)t}\|_{M_q^q}\,\|u_0\|_{L^p} \\
& \qquad \leq C\,t^{-\frac{n}\sigma\left(\frac1p-\frac1q\right)}\,e^{-ct}\,\|u_0\|_{L^p},
\end{split} \]
where we used~\eqref{eq:osclow} together with
\[ \|(1-\chi)\,e^{-a(\xi)t}\|_{M_q^q}= e^{-ct}\,\|(1-\chi)\,e^{-(a(\xi)-c)t}\|_{M_q^q} \leq C\,e^{-ct}. \]
This latter is a consequence of~\eqref{eq:oscillationshigh} for sufficiently small~$c$ so that
\begin{equation}\label{eq:aexp}
a(\xi) \geq a_1\xii^{\theta_1} \geq \frac{a_1}2\xii^{\theta_1} + c, \quad \xi\in\supp (1-\chi),
\end{equation}
%
and of Lemma~\ref{lem:hexp3} with $\phi=1-\chi$, together with Lemma~\ref{lem:exp}. %
On the other hand, if $1<p\leq q<\infty$, then
\[ \begin{split}
& \big\|\mathscr{F}^{-1}\big((1-\chi_t) e^{\pm iw(\xi)t}\,(1-\chi)\,e^{-a(\xi)t}\,\hat u_0\big)\big\|_{L^q}\\
& \qquad \leq \|(1-\chi_t)\xii^{-\beta\sigma} e^{\pm iw(\xi)t}\|_{M_p^q}\,\|(1-\chi)\,\xii^{(\beta-s)\sigma}\,e^{-a(\xi)t}\|_{M_q^q}\\
& \qquad \qquad \times\|\xii^{s\sigma}\<\xi\>^{-s\sigma}\|_{M_q^q}\|u_0\|_{H^{s\sigma,p}} \\
& \qquad \leq C\,t^{\beta-\frac{n}\sigma\left(\frac1p-\frac1q\right)}\,t^{-(\beta-s)\frac{\sigma}{\theta_1}}\,e^{-ct}\,\|u_0\|_{H^{s\sigma,p}},
\end{split} \]
where we used~\eqref{eq:M11reg} and~\eqref{eq:oschi} with $\beta=d(p,q)$ and $s$ such that
\[ d(p,q)-(d(p,q)-s)\frac{\sigma}{\theta_1}\leq0, \]
together with
\[ \begin{split}
\|(1-\chi)\,\xii^{(\beta-s)\sigma}e^{-a(\xi)t}\|_{M_q^q}
    & = e^{-ct}\,\|(1-\chi)\,\xii^{(\beta-s)\sigma}\,e^{-(a(\xi)-c)t}\|_{M_q^q} \\
    & \leq C\,t^{-(\beta-s)\frac\sigma{\theta_1}}\,e^{-ct}.
\end{split} \]
This latter is a consequence of~\eqref{eq:aexp}, as in the previous step, together with Lemma~\ref{lem:exp2} with $b=(1-\chi)\,\xii^{(\beta-s)\sigma}$ (we notice that $\beta>s$, due to $\theta_1>0$). If $p=1$, $q\in[1,\infty]$, then
\[ \begin{split}
& \big\|(1-\chi_t) e^{\pm iw(\xi)t}\,(1-\chi)\,\xii^{s\sigma}\,e^{-a(\xi)t}\big\|_{M_1^q}\\
& \qquad \leq \|(1-\chi_t)\xii^{-\beta\sigma} e^{\pm iw(\xi)t}\|_{M(H^1,L^q)}\,\|(1-\chi)\,\xii^{(\beta-s)\sigma}\,e^{-a(\xi)t}\|_{M(L^1,H^1)}\\
& \qquad \leq C\,t^{\beta-\frac{n}\sigma\left(1-\frac1q\right)}\,t^{-(\beta-s)\frac{\sigma}{\theta_1}}\,e^{-ct},
\end{split} \]
where we used~\eqref{eq:oschiH1} with $\beta=d(1,q)$ and Corollary~\ref{cor:H1} (we notice that $\beta-s>0$). By duality, the same estimate holds for $q=\infty$ and $p\in[1,\infty]$.


Therefore,
\[ \|(\id-T\chi)u(t,\cdot)\|_{L^q} \leq C\,t^{-\frac{n}\sigma\left(\frac1p-\frac1q\right)+d(p,q)-(d(p,q)-s)\frac\sigma{\theta_1}}\,e^{-ct}\,\|u_0\|_{H^{s,p}}. \]
%
%
%
If $\theta_1\geq\sigma$, we just fix $s=0$ and we get
\[  \|(\id-T\chi)u(t,\cdot)\|_{L^q} \leq C\,t^{-\frac{n}\sigma\left(\frac1p-\frac1q\right)}\,e^{-ct}\,\|u_0\|_{L^p}, \]
where we used~\eqref{eq:oschi} or~\eqref{eq:oschiH1} with $\beta=d(p,q)$ together 
with Lemma~\ref{lem:hexp3} with $\phi=1-\chi$, and with Lemma~\ref{lem:exp}. %
%
%
%
%
This concludes the proof.
\end{proof}


\section{Proof of Theorems~\ref{thm:Wavedec} and~\ref{thm:Wavereg}}
\label{sec:waveproofs}

Before giving the proof of Theorems~\ref{thm:Wavedec} and~\ref{thm:Wavereg}, the following localization argument is helpful.

\subsection{Localization of the fundamental solution}\label{sec:localization}

Let $\hat u(t,\xi)=\hat K(t,\xi)\hat u_1$ be as in \eqref{eq:wavesol}-\eqref{eq:wavesolOD}. 
For any $\chi,\chi'\in \mathcal C_c^\infty(\R^n)$, with $\chi=1$ in a neighborhood of the origin,
\[ \|(1-\chi)\,\chi' \hat K(t,\xi)\|_{M_p^q} \leq e^{-ct}, \]
for some $c>0$, since $(1-\chi)\,\chi' \hat K(t,\xi)\in \mathcal C_c^{n+1}(\R^n)$ and $a(\xi)\geq c>0$ for any $\xi\in \supp (1-\chi)\,\chi'$. Therefore, we may assume without restriction that the support of $\chi$ is contained in a sufficiently small neighborhood of the origin, or that $\chi=1$ in a sufficiently large compact. This allows us to use the expression in~\eqref{eq:wavesol} (``damped oscillations regime''), that is,
\begin{equation}\label{eq:KDO}
\hat K(t,\xi)= e^{-t\frac{a(\xi)}2}\,t\,\sinc(tw(\xi)\,\sqrt{1-\tilde a(\xi)}),\quad \text{where} \ \tilde a(\xi)=\frac{a(\xi)^2}{4w(\xi)^2},
\end{equation}
with no loss of generalization, due to $a(\xi)=\textit{o}(w(\xi))$, so that $\tilde a(\xi)=\textit{o}(1)$. Indeed, this latter property is guaranteed by the fact that $w(\xi)>0$ is homogeneous (see~\eqref{eq:poshom}), and by the assumptions $\theta_0>\sigma$ in Theorem~\ref{thm:Wavedec} and $\theta_1<\sigma$ in Theorem~\ref{thm:Wavereg}. In the cases $\theta_0<\sigma$ or $\theta_1>\sigma$, the equation would be in the ``overdamping regime'', while in the limit cases $\theta_0=\sigma$ or $\theta_1=\sigma$, the regime would possibly be mixed. 


Let us define
\[ g(\xi)=\sqrt{1-\tilde a(\xi)}-1, \quad g(\xi)=-\frac{\tilde a(\xi)}2\int_0^1\frac1{\sqrt{1-\tilde a(\xi)\tau}}\,d\tau.\]
We notice that
\[ |\partial_{\xi}^{\gamma} g(\xi)| \leq \begin{cases}
C |\xi|^{2(\theta_0-\sigma)-|\gamma|} & \text{if $\xi\in \supp\chi$,}\\
C |\xi|^{-2(\sigma-\theta_1)-|\gamma|} & \text{if $\xi\in\supp (1-\chi)$,}
\end{cases} \]
for $|\gamma|\leq n+1$. By using Taylor's formula w.r.t. $g(\xi)$, we have 
\begin{equation}\label{eq:Taylor}
\hat K(t,\xi)=\sum_{\ell=0}^N \frac1{\ell!} \left( \partial_t^{\ell} \frac{\sin(tw(\xi))}{w(\xi)}\right) \frac{e^{-ta(\xi)/2}}{1+g(\xi)} (tg(\xi))^\ell + e^{-ta(\xi)/2}R_N(t, \xi),
\end{equation}
where
\[\begin{split}
R_N(t, \xi) = \frac1{2iw(\xi) (1+g(\xi))\,N!}\int_0^1 (1-s)^N & \Big[e^{itw(\xi)(1+sg(\xi))}(itw(\xi) g(\xi))^{N+1}\\
&- e^{-itw(\xi)(1+sg(\xi))}(-itw(\xi) g(\xi))^{N+1}\Big]\,ds.
\end{split}\]
Indeed, this follows from
\[\begin{split}
e^{\pm itw(\xi) g(\xi)} & = 1 \pm itw(\xi)g(\xi)) - \frac12 (tw(\xi)g(\xi))^2+\ldots +\frac{1}{N!} \,(\pm itw(\xi)g(\xi))^N\\
&\qquad \qquad + \frac{1}{N!} (\pm itw(\xi)g(\xi))^{N+1} \,\int_0^t ( 1-s)^N e^{\pm istw(\xi)g(\xi)}\, ds,
\end{split}\]
so that
\[ \begin{split}
t\,\sinc(tw(\xi)(1+g(\xi)))
& = t\,\sum_{\ell=0}^N\partial_t^{\ell}\left(\frac{e^{itw(\xi)}-e^{-itw(\xi)}}{2i\ell!}\right)\frac{(tg(\xi))^{\ell}}{tw(\xi)(1+g(\xi))} + R_N(t, \xi) \\
& = \sum_{\ell=0}^N\frac1{\ell!}\Big(\partial_t^{\ell}\sin(tw(\xi))\Big)\frac{(tw(\xi))^{\ell-1}}{1+g(\xi)}\,(g(\xi))^{\ell}   +R_N(t, \xi).
\end{split} \]
The leading term in~\eqref{eq:Taylor} is obtained for $\ell=0$, i.e., it is
\[ t\,\sinc(tw(\xi))\,\frac{e^{-ta(\xi)/2}}{1+g(\xi)}. \]

\subsection{Proof of Theorem~\ref{thm:Wavedec}}\label{sec:proofthd}

Due to the compact support of $\chi$, it is easy to see that there exists $C>0$ such that
\[ \|\chi\hat K(t,\cdot)\|_{M_p^q} \leq C, \quad t\in[0,1], \]
for any $1\leq p\leq q\leq\infty$, so that we may assume $t\geq1$ in the following. For instance, it is sufficient to apply Lemma 8.5 in~\cite{DAE21}, as in the proof of Theorem~\ref{thm:Schdec}.


We first prove that the leading term
\[ \chi\,t\,\sinc(tw(\xi))\,\frac{e^{-ta(\xi)/2}}{1+g(\xi)} \]
verifies the desired estimate. We may follow the proof of Theorem~\ref{thm:Schdec}. On the one hand,
\[ t\,\|\chi_t\,\sinc(tw(\xi))\|_{M_p^q}\,\Big\|\chi\,\frac{e^{-ta(\xi)/2}}{1+g(\xi)}\Big\|_{M_q^q} \leq C\,t^{1-\frac{n}\sigma\left(\frac1p-\frac1q\right)}, \]
where we used~\eqref{eq:waveosclow} together with Lemma~\ref{lem:exp3} with $b=\chi/(1+g(\xi))$; indeed,
\[ \Big|\partial_\xi^\gamma \frac1{1+g(\xi)}\Big| = \Big|\partial_\xi^\gamma \frac{g(\xi)}{1+g(\xi)}\Big| \leq C\,|\xi|^{2(\theta_0-\sigma)-|\gamma|}, \qquad |\gamma|=1,\ldots,n+1. \]
On the other hand, if $1<p\leq q<\infty$, then
\begin{equation}\label{eq:maind}\begin{split}
& \|(1-\chi_t)\,\xii^{-d\sigma}\,\sin(tw(\xi))\|_{M_p^q}\,\Big\|\chi\,\xii^{d\sigma}\frac{e^{-ta(\xi)/2}}{w(\xi)(1+g(\xi))}\Big\|_{M_q^q} \\
    & \qquad \leq C\,t^{d-\frac{n}\sigma\left(\frac1p-\frac1q\right)-(d-1)\frac\sigma{\theta_0}} = C\,t^{1-\frac{n}\sigma\left(\frac1p-\frac1q\right)+(d-1)\left(1-\frac\sigma{\theta_0}\right)},
\end{split}\end{equation}
where we used~\eqref{eq:oschi} together with Lemma~\ref{lem:exp2} with
\[ b=\chi\,\xii^{d\sigma}\frac1{w(\xi)(1+g(\xi))}, \qquad \eta=(d-1)\sigma. \]
In the case $p=1$, we may use~\eqref{eq:oschiH1} together with Corollary~\ref{cor:H1} to get
\begin{equation}\label{eq:maindH1}\begin{split}
& \|(1-\chi_t)\,\xii^{-d\sigma}\,\sin(tw(\xi))\|_{M(H^1,L^q)}\,\Big\|\chi\,\xii^{d\sigma}\frac{e^{-ta(\xi)/2}}{w(\xi)(1+g(\xi))}\Big\|_{M(L^1,H^1)} \\
    & \qquad \leq C\,t^{1-\frac{n}\sigma\left(1-\frac1q\right)+(d-1)\left(1-\frac\sigma{\theta_0}\right)}.
\end{split}\end{equation}
By duality, we also obtain the estimate when $q=\infty$. %
Therefore, the proof of Theorem~\ref{thm:Wavedec} is concluded if we prove the following.
\begin{lemma}\label{lem:d}
There exists~$\delta>0$ such that
\[ \begin{split}
& \|\chi\,\left( \partial_t^{\ell} \frac{\sin(tw(\xi))}{w(\xi)}\right) \frac{e^{-ta(\xi)/2}}{1+g(\xi)} (tg(\xi))^\ell \|_{M_p^q} = t^{1-\frac{n}\sigma\left(\frac1p-\frac1q\right)+(d-1)\left(1-\frac\sigma{\theta_0}\right)-\delta}, \quad \ell\geq1, \\
& \|\chi\,e^{-ta(\xi)/2}R_N(t, \xi)\|_{M_p^q} = t^{1-\frac{n}\sigma\left(\frac1p-\frac1q\right)+(d-1)\left(1-\frac\sigma{\theta_0}\right)-\delta}, \quad \text{for sufficiently large $N$,}
\end{split} \]
for any $t\geq1$.
\end{lemma}
\begin{proof}
The reminder term $R_N$ verifies any desired polynomial decay estimate, for sufficiently large $N$. Indeed,
%
\[ |\partial_{\xi}^{\gamma}( R_N(t, \xi))| \leq C \xii^{(2\theta_0-\sigma)(N+1)-\sigma-|\gamma|}t^{N+1}.\] 
%
Due to the fact that $(2\theta_0-\sigma)(N+1)>\sigma$ for any $N\geq0$, we may apply Lemma~\ref{lem:exp2}, obtaining
\[ \|\mathscr{F}^{-1}[ \chi R_N(t, \xi)e^{-ta(\xi)/2} ] \|_{L^r} \leq C t^{-(N+1)\left(1-\frac\sigma{\theta_0}\right)-\frac{n}{\theta_0}\left(1-\frac1r\right)+\frac\sigma{\theta_0}}; \]
%
%
%
in particular, $\|\chi\,e^{-ta(\xi)/2}R_N(t, \xi)\|_{M_p^q}$ decays faster than any given polynomial, as $t\to\infty$, for sufficiently large $N$, due to $\theta_0>\sigma$.

Now we consider the expansion terms for $\ell\geq1$. On the one hand,
\[ \begin{split}
& \|\chi_t \sin^{(\ell)}(tw(\xi))\|_{M_p^q}\,t^\ell\,\Big\| \chi\, \frac{e^{-ta(\xi)/2}}{w(\xi)(1+g(\xi))} (w(\xi)g(\xi))^\ell \Big\|_{M_q^q} \\
& \qquad \leq C\,t^{\ell-\frac{n}\sigma\left(\frac1p-\frac1q\right)-2\ell+\frac{(\ell+1)\sigma}{\theta_0}} = t^{1-\frac{n}\sigma\left(\frac1p-\frac1q\right)-(\ell+1)\left(1-\frac{\sigma}{\theta_0}\right)},
\end{split} \]
where we used~\eqref{eq:waveosclow} together with Lemma~\ref{lem:exp2} with
\[ b(\xi) = \chi\, \frac{(w(\xi)g(\xi))^\ell}{w(\xi)(1+g(\xi))}, \qquad \eta=2\ell\theta_0 - (\ell+1)\sigma.  \]
In particular, $\eta>0$ due to $\ell\geq1$ and $\theta_0>\sigma$. On the other hand, if $1<p\leq q<\infty$,
\[ \begin{split}
& \|(1-\chi_t)\,\xii^{-d\sigma}\,\sin^{(\ell)}(tw(\xi))\|_{M_p^q}\,t^\ell\,\Big\| \chi\,\xii^{d\sigma}\,\frac{e^{-ta(\xi)/2}}{w(\xi)(1+g(\xi))} (w(\xi)g(\xi))^\ell \Big\|_{M_q^q} \\
& \qquad \leq C\,t^{\ell+d-\frac{n}\sigma\left(\frac1p-\frac1q\right)-2\ell+(\ell+1-d)\frac{\sigma}{\theta_0}} = C\,t^{1-\frac{n}\sigma\left(\frac1p-\frac1q\right)+(d-\ell-1)\left(1-\frac{\sigma}{\theta_0}\right)},
\end{split} \]
where we used~\eqref{eq:oschi}, together with Lemma~\ref{lem:exp2} with
\[ b(\xi) = \chi\,\xii^{d\sigma}\,\frac{(w(\xi)g(\xi))^\ell}{w(\xi)(1+g(\xi))}, \qquad \eta=2\ell\theta_0 + (d-\ell-1)\sigma.  \]
We get the desired result for any $\ell\geq1$. When $p=1$ or $q=\infty$, we proceed in the usual way, replacing~\eqref{eq:oschi} and Lemma~\ref{lem:exp2} with~\eqref{eq:oschiH1} and Corollary~\ref{cor:H1}. This concludes the proof.
\end{proof}
\begin{Rem}\label{rem:sharpdecay}
Thanks to Lemma~\ref{lem:d}, Theorem~\ref{thm:Wavedec} is sharp, as far as~\eqref{eq:maind} is sharp. We explicitly prove the optimality of this latter in the radial case in Lemma~\ref{lem:crucialwave}.
\end{Rem}

\subsection{Proof of Theorem~\ref{thm:Wavereg}}\label{sec:proofthr}

We may follow the proof of Theorem~\ref{thm:Wavedec}, with some modifications. %
%
%
A crucial difference is that now
\[ e^{-t\frac{a(\xi)}2} = e^{-c_1t}\,e^{-t\frac{a(\xi)-2c_1}2}, \]
where $a(\xi)-2c_1$ still verifies the same assumption~\eqref{eq:oscillationshigh} on $a(\xi)$ for sufficiently small $c_1$, due to
\[ a(\xi)-2c_1 \geq a_1\xii^{\theta_1}-2c_1 \geq \frac{a_1}2\,\xii^{\theta_1},\quad \text{for sufficiently large $\xii$.} \]
Since we produced exponential decay in time, we are only interested in a regularity result and in controlling the possible singular behavior of the estimate as $t\to0$.

We first prove that the leading term
\[ (1-\chi)\,t\,\sinc(tw(\xi))\,\frac{e^{-ta(\xi)/2}}{1+g(\xi)} \]
verifies the desired estimate. We may follow the proof of Theorem~\ref{thm:Schreg}. First, let $\theta_1\in(0,\sigma)$. On the one hand,
\[ t\,\|\chi_t\,\sinc(tw(\xi))\|_{M_p^q}\,\Big\|(1-\chi)\,\frac{e^{-ta(\xi)/2}}{1+g(\xi)}\Big\|_{M_q^q} \leq C\,t^{1-\frac{n}\sigma\left(\frac1p-\frac1q\right)}, \]
thanks to~\eqref{eq:waveosclow} 
%
%
and of Lemma~\ref{lem:hexp3} with
\[ \phi(\xi)=(1-\chi)\,\frac1{1+g(\xi)}, \]
together with Lemma~\ref{lem:exp}. %
On the other hand, if $1<p\leq q<\infty$, then
\begin{equation}\label{eq:mainr}\begin{split}
& \|(1-\chi_t)\,\xii^{-d\sigma}\,\sin(tw(\xi))\|_{M_p^q}\,\Big\|(1-\chi)\,\xii^{(d-s)\sigma}\frac{e^{-ta(\xi)/2}}{w(\xi)(1+g(\xi))}\Big\|_{M_q^q} \\
    & \qquad \leq C\,t^{d-\frac{n}\sigma\left(\frac1p-\frac1q\right)-(d-s-1)\frac\sigma{\theta_1}} = C\,t^{1-\frac{n}\sigma\left(\frac1p-\frac1q\right)+(d-1)-(d-1-s)\frac\sigma{\theta_1}},
\end{split}\end{equation}
thanks to~\eqref{eq:oschi} and Lemma~\ref{lem:exp2} with
\[ b(\xi)=(1-\chi)\,\xii^{(d-s)\sigma}\frac{1}{w(\xi)(1+g(\xi))}, \quad \eta=(d-s-1)\sigma. \]
We notice that $d-1>s$, due to $\theta_1>0$. As in the proof of Theorem~\ref{thm:Schreg},
\[ \|\xii^{s\sigma}\,\<\xi\>^{-s\sigma}\|_{M_q^q}\,\|\mathscr{F}^{-1}(\<\xi\>^{s\sigma}\hat u_1)\|_{L^p}\leq C\,\|u_1\|_{H^{s\sigma,p}}. \]
%

If $p=1$ and $q\in[1,\infty]$, we may use~\eqref{eq:oschiH1} together with Corollary~\ref{cor:H1} to get
\begin{equation}\label{eq:mainrH1}\begin{split}
& \|(1-\chi_t)\,\xii^{-d\sigma}\,\sin(tw(\xi))\|_{M(H^1,L^q)}\,\Big\|(1-\chi)\,\xii^{(d-s)\sigma}\frac{e^{-ta(\xi)/2}}{w(\xi)(1+g(\xi))}\Big\|_{M(L^1,H^1)} \\
    & \qquad \leq C\,t^{1-\frac{n}\sigma\left(1-\frac1q\right)+(d-1)-(d-1-s)\frac\sigma{\theta_1}}.
\end{split}\end{equation}
By duality, we obtain the estimate when $q=\infty$. %

Since producing the exponential decay for the other terms of the Taylor expansion is obvious, the proof of Theorem~\ref{thm:Wavereg} is concluded if we prove the following.
\begin{lemma}\label{lem:r}
There exists $\delta>0$ such that
\[ \begin{split}
& \|(1-\chi)\,\left( \partial_t^{\ell} \frac{\sin(tw(\xi))}{w(\xi)}\right) \frac{e^{-ta(\xi)/2}}{1+g(\xi)} (tg(\xi))^\ell \|_{M_p^q} \\
& \qquad = t^{1-\frac{n}\sigma\left(\frac1p-\frac1q\right)+(d-1)-(d-1-s)\frac\sigma{\theta_1}+\delta}, \quad \ell\geq1, \\
& \|(1-\chi)\,e^{-ta(\xi)/2}R_N(t, \xi)\|_{M_p^q} \\
& \qquad = t^{1-\frac{n}\sigma\left(\frac1p-\frac1q\right)+(d-1)-(d-1-s)\frac\sigma{\theta_1}+\delta}, \quad \text{for sufficiently large $N$,}
\end{split} \]
for any $t\in(0,1]$.
\end{lemma}
%
%
\begin{proof}
Since we are interested in $t\to0$, we restrict to $t\in(0,1]$. 
We will use~\eqref{eq:Taylor} and we proceed as we did in the proof of Theorem~\ref{thm:Wavedec}. We have
\[ |\partial_{\xi}^{\gamma}R_N(t, \xi)|\lesssim \xii^{(2\theta_1-\sigma)(N+1)-\sigma-|\gamma|}t^{N+1}. \]
If $\theta_1\in(\sigma/2,\sigma)$, we can apply Lemma~\ref{lem:exp2} for sufficiently large $N$, so that
\[\begin{split}
\|(1-\chi)\,R_N(t, \xi)e^{-ta(\xi)/2}\|_{M_p^q}
    & \leq\|\mathscr{F}^{-1}( (1-\chi)\,R_N(t, \xi)e^{-ta(\xi)/2})\|_{L^r}\\
    & \lesssim t^{(N+1)\left(\frac\sigma{\theta_1}-1\right)+\frac1{\theta_1}\left(\sigma-n\left(\frac1p-\frac1q\right)\right)}.
\end{split} \]
For sufficiently large $N$, this gives the desired estimate. If $\theta_1\in(0,\sigma/2]$, 
%
we may weaken the estimate for the derivative to
\[ |\partial_{\xi}^{\gamma}R_N(t, \xi)|\lesssim \xii^{\eta-|\gamma|}t^{N+1}, \]
for some $\eta>0$. Then we apply Lemma~\ref{lem:exp2}, obtaining
\[ \|(1-\chi)\,R_N(t, \xi)e^{-ta(\xi)/2}\|_{M_p^q} \lesssim t^{N+1+\frac1{\theta_1}\left(-\eta-n\left(\frac1p-\frac1q\right)\right)}. \]
%
For sufficiently large $N$, this gives the desired estimate.

Now we consider the expansion terms for $\ell\geq1$. On the one hand, we estimate
\[ \|\chi_t \sin^{(\ell)}(tw(\xi))\|_{M_p^q} = t^{-\frac{n}\sigma\left(\frac1p-\frac1q\right)},\]
using~\eqref{eq:waveosclow}, while we estimate
\[ t^\ell\,\Big\| (1-\chi)\, \frac{e^{-ta(\xi)/2}}{w(\xi)(1+g(\xi))} (w(\xi)g(\xi))^\ell \Big\|_{M_q^q} \leq C t^{1+(\ell+1)\left(\frac{\sigma}{\theta_1}-1\right)}, \]
by using Lemma~\ref{lem:exp2} with
\[ b(\xi) = (1-\chi)\, \frac{(w(\xi)g(\xi))^\ell}{w(\xi)(1+g(\xi))},  \]
and letting $\eta= 2\ell\theta_1- (\ell+1)\sigma$ if this latter is positive. Otherwise, we fix $\eta\in(0,\ell\theta_1)$, and we use Lemma~\ref{lem:exp2} to estimate
\[ t^\ell\,\Big\| (1-\chi)\, \frac{e^{-ta(\xi)/2}}{w(\xi)(1+g(\xi))} (w(\xi)g(\xi))^\ell \Big\|_{M_q^q} \leq C t^{\ell-\frac\eta{\theta_1}}. \]
On the other hand, if $1<p\leq q<\infty$, then
\[ \begin{split}
& \|(1-\chi_t)\,\xii^{-d\sigma}\,\sin^{(\ell)}(tw(\xi))\|_{M_p^q}\,t^\ell\,\Big\| (1-\chi)\,\xii^{(d-s)\sigma}\,\frac{e^{-ta(\xi)/2}}{w(\xi)(1+g(\xi))} (w(\xi)g(\xi))^\ell \Big\|_{M_q^q} \\
& \qquad \leq C t^{-\frac{n}\sigma\left(\frac1p-\frac1q\right)+d-(d-s-1)\frac{\sigma}{\theta_1} + \ell\left(\frac{\sigma}{\theta_1}-1\right)},
\end{split} \]
where we used~\eqref{eq:oschi}, together with Lemma~\ref{lem:exp2} with
\[ b(\xi) = \chi\,\xii^{(d-s)\sigma}\,\frac{(w(\xi)g(\xi))^\ell}{w(\xi)(1+g(\xi))},  \]
and letting $\eta=2\ell\theta_1 + (d-s-\ell-1)\sigma$, when $\eta$ is positive, so we get the desired result. 
%
%
Otherwise, we proceed as before, choosing $\eta>0$, sufficiently small. When $p=1$ or $q=\infty$, we proceed in the usual way, replacing~\eqref{eq:oschi} and Lemma~\ref{lem:exp2} with~\eqref{eq:oschiH1} and Corollary~\ref{cor:H1}. %
%
This concludes the proof.
\end{proof}
\begin{Rem}\label{rem:sharpreg}
Thanks to Lemma~\ref{lem:r}, Theorem~\ref{thm:Wavereg} is sharp, as far as~\eqref{eq:mainr} is sharp. We explicitly prove the optimality of this latter in the radial case in Lemma~\ref{lem:crucialwave}.
\end{Rem}


\subsection{Proofs with $L^2$ regularity}

In order to prove Proposition~\ref{prop:Wavedec2}, it is sufficient to modify the proof of Theorem~\ref{thm:Wavedec} in~\textsection\ref{sec:proofthd} for the leading term
\[ \chi\,t\,\sinc(tw(\xi))\,\frac{e^{-ta(\xi)/2}}{1+g(\xi)}, \]
for $t\geq1$, as follows. Thanks to $M_1^2=L^2$ and $n=2\sigma$, using  the change of variable~$t^{\frac1{\theta_0}}\xi\mapsto\xi$, we now estimate
\[ \begin{split}
t\,\Big\|(1-\chi_t)\,\chi\,\sinc(tw(\xi))\,\frac{e^{-ta(\xi)/2}}{1+g(\xi)}\Big\|_{L^2}
    & \leq C\Big(\int_{\xii\geq Mt^{-\frac1\sigma}} \xii^{-n}\,e^{-2a_0t\xii^{\theta_0}}\,d\xi\Big)^{\frac12} \\
    & = C \Big(\int_{\xii\geq Mt^{\frac1{\theta_0}-\frac1\sigma}} \xii^{-n}\,e^{-2a_0\xii^{\theta_0}}\,d\xi\Big)^{\frac12} \\
    & \leq C_1\,(1+\log t)^{\frac12}.
\end{split} \]
Similarly, to prove Proposition~\ref{prop:Wavereg2}, it is sufficient to modify the proof of Theorem~\ref{thm:Wavereg} in~\textsection\ref{sec:proofthr} for the leading term
\[ (1-\chi)\,t\,\sinc(tw(\xi))\,\frac{e^{-ta(\xi)/2}}{1+g(\xi)}\,, \]
for $t\in(0,1]$ as follows. Thanks to $M_1^2=L^2$ and $n=2\sigma$, using  the change of variable~$t^{\frac1{\theta_1}}\xi\mapsto\xi$, we now estimate 
\[ \begin{split}
t\,\Big\|(1-\chi_t)\,(1-\chi)\,\sinc(tw(\xi))\,\frac{e^{-ta(\xi)/2}}{1+g(\xi)}\Big\|_{L^2}
    & \leq \Big(\int_{\xii\geq Mt^{-\frac1\sigma}} \xii^{-n}\,e^{-2a_1t\xii^{\theta_1}}\,d\xi\Big)^{\frac12} \\
    & = C \Big(\int_{\xii\geq Mt^{\frac1{\theta_1}-\frac1\sigma}} \xii^{-n}\,e^{-2a_1\xii^{\theta_1}}\,d\xi\Big)^{\frac12} \\
    & \leq C\,(1-\log t)^{\frac12}.
\end{split} \]
We stress that the different sign in front of the $\log t$ term depends on the sign of
\[ -\log t^{\frac1{\theta_j}-\frac1\sigma}= \left(\frac1\sigma-\frac1{\theta_j}\right) \log t, \quad j=0,1. \]


\section{Optimality of the estimates}\label{sec:optimality}

The optimality of the estimates is expected from the proofs, since it is related to the different scaling of $e^{\pm iw(\xi)}$ or $\sinc w(\xi)$, and of $e^{-a(\xi)}$. However, we can provide a direct proof of optimality, at least in the radial case. That is, we assume that $w(\xi)=c\xii^\sigma$ for some $c>0$ (we may assume $c=1$ with no loss of generality) and that $a(\xi)=a_0\xii^{\theta_0}$, for some $\theta_0>\sigma$, or $a(\xi)=a_1\xii^{\theta_1}$, for some $\theta_1\in(0,\sigma)$.

We may indeed relax the radial assumption for $a(\xi)$. For instance, it is sufficient to assume that there exists $r\in\mathcal C^{n+1}$ such that
\[ a(\xi) = a_0\xii^{\theta_0}\,(1+r(\xi))\,,\qquad |\partial_\gamma r(\xi)|\leq C\,\xii^{\varepsilon-|\gamma|},\]
in a neighborhood of the origin, or
\[ a(\xi) = a_1\xii^{\theta_1}\,(1+r(\xi))\,,\qquad |\partial_\gamma r(\xi)|\leq C\,\xii^{-\varepsilon-|\gamma|},\]
out of some compact, for some $\varepsilon>0$. Indeed, it is easy to prove that the solution to~\eqref{eq:CPSch} and to~\eqref{eq:CPwave} with $a(\xi)=a_0\xii^{\theta_0}\,r(\xi)$ verifies better estimates than the ones provided by Theorems~\ref{thm:Schdec} and~\ref{thm:Wavedec} for $a(\xi)=a_0\xii^{\theta_0}$, and that the solution to~\eqref{eq:CPSch} and to~\eqref{eq:CPwave} with $a(\xi)=a_1\xii^{\theta_1}\,r(\xi)$ verifies better estimates than the ones provided by Theorems~\ref{thm:Schreg} and~\ref{thm:Wavereg} for $a(\xi)=a_1\xii^{\theta_1}$.

As expected, the crucial point of the optimality is based on the optimality of the estimates
\begin{align}
\label{eq:optdec}
\|(1-\chi_t)\,\xii^{-s\sigma}\,e^{\pm i\xii^\sigma t-a_0\xii^{\theta_0}t} \|_{M_p^q}
    & \leq t^{-\frac{n}\sigma\left(\frac1p-\frac1q\right)+s+(d-s)\left(1-\frac\sigma{\theta_0}\right)}, \qquad t\geq1, \\
\label{eq:optreg}
\|(1-\chi_t)\,\xii^{-s\sigma}\,e^{\pm i\xii^\sigma t-a_1\xii^{\theta_1}t} \|_{M_p^q}
    & \leq t^{-\frac{n}\sigma\left(\frac1p-\frac1q\right)+s-(d-s)\left(\frac\sigma{\theta_1}-1\right)}, \qquad t\in(0,1],
\end{align}
for $s\in[0,d]$, where $d$ is as in~\eqref{eq:dpqdef}. Indeed, the optimality of~\eqref{eq:Schdec} and~\eqref{eq:Wavedec} follows from~\eqref{eq:optdec} (for $s=0$ and $s=1$, respectively), while the optimality of~\eqref{eq:Schreg} and~\eqref{eq:Wavereg} follows from~\eqref{eq:optreg} (replacing $s$ by $s+1$ for the latter).

By the change of variable $t\xii^\sigma\mapsto\xii^\sigma$, the previous estimates reduce to
\[ \begin{split}
\|(1-\chi)\,\xii^{-s\sigma}\,e^{\pm i\xii^\sigma-a_0\xii^{\theta_0}t^{-\frac{\theta_0}\sigma+1}} \|_{M_p^q}
    & \leq t^{(d-s)\left(1-\frac\sigma{\theta_0}\right)}, \qquad t\geq1, \\
\|(1-\chi)\,\xii^{-s\sigma}\,e^{\pm i\xii^\sigma-a_1\xii^{\theta_1}t^{1-\frac{\theta_1}\sigma}} \|_{M_p^q}
    & \leq t^{(d-s)\left(1-\frac\sigma{\theta_1}\right)}, \qquad t\in(0,1],
\end{split} \]
that is, setting $\tau=a_j^{\frac1\theta}\,t^{\frac1{\theta}-\frac1\sigma}$ where we put $\theta=\theta_0$ in the first case and $\theta=\theta_1$ in the second one, we shall prove the optimality of the singular estimate
\[ \|(1-\chi)\,\xii^{-s\sigma}\,e^{\pm i\xii^\sigma-(\xii\tau)^\theta} \|_{M_p^q} \leq \tau^{-(d-s)\sigma}, \qquad \tau\in(0,1]. \]
\begin{lemma}\label{lem:crucialwave}
Let $n\geq1$, $\sigma>0$, $\theta>0$, $1\leq p \leq q\leq \infty$ and $d$ as in~\eqref{eq:dpqdef}. Assume that $s\in[0,d]$. Then there exist $C>0$ such that
\begin{equation}\label{eq:crucialwave}
\| (1-\chi)\,\xii^{-s\sigma}\,e^{\pm i\xii^\sigma-(\xii\tau)^\theta} \|_{M_p^q}\geq C\,\tau^{-(d-s)\sigma}\,,\quad \tau\in(0,1].
\end{equation}
\end{lemma}
\begin{proof}
By duality, we assume without restriction that $p\leq q'$.

The case $\sigma\neq1$ is proved in~\cite[Proposition 7.2]{DAE21}, so we only need to prove the statement for $\sigma=1$. We fix $g\in \mathcal C_c^\infty$, radial, and we define $g_\tau(x)= \tau^{-\frac{n}p}\, g(\tau^{-1} x)$, so that $\|g_\tau\|_{L^p}=\|g\|_{L^p}$. Then $\widehat{g_\tau} (\xi) = \tau^{\frac{n}{p'}}\,\hat g (\tau\xi)$. Using the representation formula of the Fourier transform for radial functions, we have (abusing notation, we write $\chi(\rho)$)
\[ \begin{split}
& \mathscr{F}^{-1}((1-\chi(\xi))\,\xii^{-s}\,e^{\pm i\xii-(\tau\xii)^\theta}g_\tau) (|x|) \\
& \qquad = \int_0^\infty (1-\chi(\rho))\, \rho^{-s}\,e^{\pm i\rho-(\tau\rho)^\theta}\widehat{g_\tau}(\rho)\,J_{\frac{n-2}2}(|x|\rho)\,|x|^{1-\frac{n}2}\,\rho^{\frac{n}2}\,d\rho,
\end{split} \]
where $J_{\frac{n-2}2}$ denotes the Bessel function of first kind. Now let $|x|\in I=[1-\delta\tau,1+\delta\tau]$, for some $\delta\ll \tau^{-1}$ which we will fix later. Due to the fact that $1-\chi(\rho)$ vanishes for $\rho\leq1$ and $|x|\approx 1$, inside the integral we may use the asymptotic expansion of the Bessel function
\[ J_{\nu}(z) =c_{\nu} z^{-\frac12}\,\cos\left(z-\nu\frac{\pi}{2}-\frac{\pi}{4}\right) +\textit{O}\,(|z|^{-3/2}), \quad |z|\rightarrow \infty. \]
The leading term in the asymptotic expansion above may be estimated if we estimate
\[ \begin{split}
|x|^{-\frac{n-1}2}\tau^{\frac{n}{p'}}&\int_0^\infty (1-\chi(\rho))\,e^{i\rho(1\pm |x|)}\,\hat g (\tau\rho)\,\,\rho^{\frac{n-1}2-s}e^{-(\tau\rho)^\theta}\,d\rho\\
&= |x|^{-\frac{n-1}2}\tau^{\left(\frac{n}{p'}-\frac{n+1}2+s\right)}
\int_0^\infty (1-\chi(s\tau^{-1}))\,e^{is\tau^{-1}(1\pm |x|)}\,\hat g(s)\,\,s^{\frac{n-1}2-s}e^{-s^{\theta}}\,ds\\
&=|x|^{-\frac{n-1}2}\tau^{\left(\frac{n}{p'}-\frac{n+1}2+s\right)}\,\big((2\pi)^n\,h(\tau^{-1}(1\pm |x|))-J(\tau)\big),
\end{split}
\]
where
\[\hat h(t)=\hat g(t)\,\,t^{\frac{n-1}2-s}e^{-ct^{\theta}},\quad J(\tau)=\int_0^\infty \chi(t\tau^{-1})\,e^{it\tau^{-1}(1\pm |x|)}\,\hat g(t)\,\,t^{\frac{n-1}2-s}e^{-t^{\theta}}\,dt.\]
In the last equality, we may assume without restriction that $\hat g(t)\,\,t^{\frac{n-1}2-s}$ is integrable near $t=0$, so that
\[ |J(\tau)|\leq \int_0^{2\tau} |\hat g(t)|\,\,t^{\frac{n-1}2-s}e^{-t^{\theta}}\,dt \]
%
vanishes as $\tau\to0$. By Riemann--Lebesgue theorem, $h(\tau^{-1}(1+ |x|))$ vanishes as $\tau\to0$, due to $\hat h \in L^1$, for any given $x$. Choosing $g$ such that
\[ h(0)=(2\pi)^{-n}\,\int_0^\infty \hat g(t)\,\,t^{\frac{n-1}2-s}e^{-t^{\theta}}\,dt\neq 0,\]
there exists $\delta>0$ such that $|h(\tau^{-1}(1-|x|))|\geq |h(0)|/2>0$ for $|x|\in I=[1-\delta \tau,  1+\delta \tau]$. Therefore
\[
\begin{split}
& \|\mathscr{F}^{-1}((1-\chi(\xi))\,\xii^{-s}\,e^{\pm i\xii-(\tau\xii)^\theta}g_\tau) \|_{L^q} \\
& \qquad \qquad \geq\|\mathscr{F}^{-1}((1-\chi(\xi))\,\xii^{-s}\,e^{\pm i\xii-(\tau\xii)^\theta}g_\tau) \|_{L^q(I)} \\
& \qquad \qquad \geq C\,\tau^{\left(\frac{n}{p'}-\frac{n+1}2+s\right)}\,\Big(\int_{|x|\in I} 1\,dx\Big)^{\frac1q} \approx \tau^{\left(s+\frac{n}{p'}-\frac{n+1}2+\frac1q\right)}= \tau^{s-d},
\end{split}\]
with $d$ as in~\eqref{eq:dpqdef}. This concludes the proof.
\end{proof}

We may easily prove the optimality of the estimates in Propositions~\ref{prop:Wavedec2} and~\ref{prop:Wavereg2}. In this case, a logarithmic loss appears with respect to what provided by Lemma~\ref{lem:crucialwave}, but such loss is also sharp, in the sense that the estimate from below in Lemma~\ref{lem:crucialwave} may be improved.
\begin{lemma}\label{lem:crucialwave2}
Let $n\geq1$ and $\sigma=n/2$, $\theta>0$. Then there exist $C>0$ and $\tau_0>0$ such that
\begin{equation}\label{eq:crucialwave2}
\| (1-\chi)\,\sinc(\xii^{\sigma})\,e^{-a(\xii\tau)^\theta} \|_{M_1^2}\geq C\,(1-\log\tau)^{\frac12}\,,\quad \tau\in(0,\tau_0].
\end{equation}
\end{lemma}
\begin{proof}
For a given $\chi$, there exists $\tau_0\leq1$ be such that $\supp\chi\subset [0,\tau_0^{-1})$. Let $\tau\in(0,\tau_0]$. 
%
%
%
%
Noticing that $e^{-2a(\xii\tau)^\theta}\geq e^{-2a}$ for $|\xi|\leq\tau^{-1}$, we may estimate
\[ \begin{split}
& \| (1-\chi)\,\xii^{-\sigma}\,\sin(\xii^\sigma)\,e^{-a(\xii\tau)^\theta} \|_{M_1^2}\\
    & \qquad = \Big(\int_{\R^n} (1-\chi)^2\,\xii^{-n}\,\sin^2(\xii^\sigma)\,e^{-2a(\xii\tau)^\theta}\,d\xi\Big)^{\frac12} \\
    & \qquad\geq e^{-a}\,\Big(\int_{\{\xi\in\R^n: \tau_0^{-1}\leq\xii\leq \tau^{-1}\}} \xii^{-n}\,\sin^2(\xii^\sigma)d\xi\Big)^{\frac12} \\
    & \qquad= e^{-a}|\S^{n-1}|\Big(\int_{\tau_0^{-1}}^{\tau^{-1}} \rho^{-1}\,\sin^2(\rho^\sigma)\,d\rho\Big)^{\frac12}\\
    & \qquad\geq C\,(\log\tau_0-\log\tau)^{\frac12} \geq C_1(1-\log\tau)^{\frac12},
\end{split} \]
for sufficiently small $\tau$ (we use that $\sin^2(\rho^\sigma)\geq 1/\sqrt{2}$ when $\rho^\sigma\in[(\pi/4)+m\pi,(3\pi/4)+m\pi]$, for $m\in\N$). This concludes the proof.
\end{proof}

\section{Proofs of Lemmas \ref{lem:exp}--\ref{lem:hexp3}}\label{sec:multproofs}

As discussed in \textsection\ref{sec:multipliers}, we will prove the multipliers Lemmas~\ref{lem:exp}, \ref{lem:exp2}, \ref{lem:exp3} 
and~\ref{lem:hexp3}, using the integration by parts method. Let $h\in \mathcal C^1(\R^n\setminus\{0\})$, be in $L^1(\R^n)$ with its derivatives (as a consequence, $h(\xi)$ vanishes as $\xii\to\infty$). Then
\[ (2\pi)^{n} \mathscr{F}^{-1}h = \int_{\R^n} e^{ix\xi} h(\xi)\,d\xi = i|x|^{-2} \sum_{j=1}^n x_j \int_{\R^n} e^{ix\xi} \partial_{\xi_j}h(\xi)\,d\xi. \]
Iterating this process, if $h\in \mathcal C^\kappa(\R^n\setminus\{0\})$, in $L^1(\R^n)$ with all its derivatives, then
\begin{equation}\label{eq:intparts}
(2\pi)^{n}\mathscr{F}^{-1}h = |x|^{-\kappa} \sum_{|\gamma|=\kappa} c_\gamma (ix/|x|)^\gamma \int_{\R^n} e^{ix\xi} \partial_{\xi}^\gamma h(\xi)\,d\xi.
\end{equation}
In the case in which $h\in \mathcal C^1(\R^n\setminus\{0\})$, is in $L^1(\R^n)$, and its derivatives are in $L^1(\R^n\setminus B(0,\delta))$, it may be useful to perform integration by parts only in $\R^n\setminus B(0,\delta)$, obtaining
\begin{equation}\label{eq:splitintparts}
\begin{split}
(2\pi)^{n} \mathscr{F}^{-1}h
     & = \int_{B(0,\delta)} e^{ix\xi} h(\xi)\,d\xi + i|x|^{-2} \sum_{j=1}^n x_j \int_{\R^n\setminus \bar B(0,\delta)} e^{ix\xi} \partial_{\xi_j}h(\xi)\,d\xi \\
     & \qquad - i|x|^{-2} \sum_{j=1}^n x_j \int_{\partial B(0,\delta)} e^{ix\xi} h(\xi)\,\nu_j\,dS(\xi).
\end{split}
\end{equation}
Employing~\eqref{eq:intparts} and~\eqref{eq:splitintparts}, we may prove the results stated in \textsection\ref{sec:multipliers}.

We first prove Lemma~\ref{lem:exp}.
\begin{proof}[Proof of Lemma~\ref{lem:exp}]
By the inversion formula of Fourier transform, it is clear that
\[ |\mathscr{F}^{-1} (e^{-ta(\xi)})|\leq (2\pi)^{-n}\,\int_{\R^n} e^{-ta(\xi)}\,d\xi  \lesssim \int_{\R^n} e^{-ct\xii^\theta}\,d\xi = Ct^{-\frac{n}\theta}, \]
where we used the change of variable $t^{\frac1\theta}\xi\mapsto\xi$. We fix $\eps\in(0,1)$ such that~$\eps\leq\theta$. For any $|\gamma|\geq1$, we may estimate
\[ |\partial_\xi^\gamma e^{-ta(\xi)}|\lesssim t\,\xii^{\theta-|\gamma|}\big(1+t\xii^\theta\big)^{|\gamma|-1}\,e^{-ct\xii^\theta}\lesssim  t^{\frac\eps\theta}\,\xii^{\eps-|\gamma|}\,e^{-\frac{c}2t\xii^\theta}. \]
In particular, $\partial_\xi^\gamma e^{-ta(\xi)}$ is in $L^1$ for $|\gamma|\leq n$, since $\eps>0$. Integrating by parts $n$ times as in~\eqref{eq:intparts} and then once as in~\eqref{eq:splitintparts} with $\delta=|x|^{-1}$; we obtain:
\[
\begin{split}
|\mathscr{F}^{-1}(e^{-ta(\xi)})|
     & \lesssim t^{\frac\eps\theta}\,|x|^{-n} \int_{B(0,|x|^{-1})} \xii^{\eps-n} d\xi \\
     & \qquad + t^{\frac\eps\theta}\,|x|^{-(n+1)} \int_{\R^n\setminus \bar B(0,|x|^{-1})} \xii^{\eps-(n+1)} \,d\xi \\
     & \qquad + t^{\frac\eps\theta}\,|x|^{-(n+1)} \int_{\partial B(0,|x|^{-1})} \xii^{\eps-n}\,dS(\xi)\\
     & \lesssim t^{\frac\eps\theta}\,|x|^{-(n+\eps)}.
\end{split}
\]
Then:
\[\begin{split}
\|\mathscr{F}^{-1} (e^{-ta(\xi)})\|_{L^r}
     &\lesssim t^{-\frac{n}\theta}\,\Big(\int_{|x|\leq t^{\frac1\theta}}1dx\Big)^{\frac1r} + t^{\frac\eps\theta} \Big(\int_{|x|\geq t^{\frac1\theta}} |x|^{-(n+\eps)r}\,dx\Big)^{\frac1r}\\
     & \lesssim t^{-\frac{n}\theta\left(1-\frac1r\right)}.
\end{split} \]
This concludes the proof.
\end{proof}
We then prove Lemma~\ref{lem:exp2}. 
This lemma is a variant of \cite[Lemma~8.5]{DAE21}, which takes into account of exponential terms.
\begin{proof}[Proof of Lemma~\ref{lem:exp2}]
By the inversion formula of Fourier transform, it is clear that
\[ |\mathscr{F}^{-1} (b(\xi)\,e^{-ta(\xi)})|\leq (2\pi)^{-n}\,\int_{\R^n} |b(\xi)|\,e^{-ta(\xi)}\,d\xi  \lesssim \int_{\R^n} \xii^\eta\,e^{-ct\xii^\theta}\,d\xi = Ct^{-\frac{n+\eta}\theta}, \]
where we used the change of variable $t^{\frac1\theta}\xi\mapsto \xi'$.

We now distinguish two cases. First, let~$\eta>0$. We fix $\eps\in(0,1)$ such that~$\eps\leq\eta$. For any $|\gamma|\geq1$, we may estimate
\[ |\partial_\xi^\gamma (b(\xi)\,e^{-ta(\xi)})|\lesssim \xii^{\eta-|\gamma|}\big(1+t\xii^\theta\big)^{|\gamma|}\,e^{-ct\xii^\theta}\lesssim  t^{-\frac{\eta-\eps}\theta}\,\xii^{\eps-|\gamma|}\,e^{-\frac{c}2t\xii^\theta}. \]
In particular, $\partial_\xi^\gamma (b(\xi)\,e^{-ta(\xi)})$ is in $L^1$ for $|\gamma|\leq n$, since $\eps>0$. We now integrate by parts $n$ times as in~\eqref{eq:intparts} and then one time as in~\eqref{eq:splitintparts} with $\delta=|x|^{-1}$; we obtain:
\[
\begin{split}
|\mathscr{F}^{-1}(b(\xi)\,e^{-ta(\xi)})|
     & \lesssim t^{-\frac{\eta-\eps}\theta}\,|x|^{-n} \int_{B(0,|x|^{-1})} \xii^{\eps-n} d\xi \\
     & \qquad + t^{-\frac{\eta-\eps}\theta}\,|x|^{-(n+1)} \int_{\R^n\setminus \bar B(0,|x|^{-1})} \xii^{\eps-(n+1)} \,d\xi \\
     & \qquad + t^{-\frac{\eta-\eps}\theta}\,|x|^{-(n+1)} \int_{\partial B(0,|x|^{-1})} \xii^{\eps-n}\,dS(\xi)\\
     & \lesssim t^{-\frac{\eta-\eps}\theta}\,|x|^{-(n+\eps)}.
\end{split}
\]
Then:
\[\begin{split}
\|\mathscr{F}^{-1} (b(\xi)\,e^{-ta(\xi)})\|_{L^r}
     &\lesssim t^{-\frac{n+\eta}\theta}\,\Big(\int_{|x|\leq t^{\frac1\theta}}1dx\Big)^{\frac1r} + t^{-\frac{\eta-\eps}\theta} \Big(\int_{|x|\geq t^{\frac1\theta}} |x|^{-(n+\eps)r}\,dx\Big)^{\frac1r}\\
     & \lesssim t^{-\frac{n}\theta\left(1-\frac1r\right)-\frac\eta\theta}.
\end{split} \]
This concludes the proof for $\eta>0$. Now let~$\eta\in(-n,0]$; for any $|\gamma|\geq1$, we may estimate
\[ |\partial_\xi^\gamma (b(\xi)\,e^{-ta(\xi)})|\lesssim \xii^{\eta-|\gamma|}\,e^{-ct\xii^\theta}. \]
Let $\kappa$ be the largest integer such that $\eta-\kappa>-n$. Assume for a moment that $\eta$ is not an integer. Since $\eta$ is not an integer, then $\eta-\kappa-1<-n$. We now integrate by parts $\kappa$ times as in~\eqref{eq:intparts} and then one time as in~\eqref{eq:splitintparts} with $\delta=|x|^{-1}$; we obtain:
\[
\begin{split}
|\mathscr{F}^{-1}(b(\xi)\,e^{-ta(\xi)})|
     & \lesssim |x|^{-\kappa} \int_{B(0,|x|^{-1})} \xii^{\eta-\kappa} d\xi \\
     & \qquad + |x|^{-(\kappa+1)} \int_{\R^n\setminus \bar B(0,|x|^{-1})} \xii^{\eta-(\kappa+1)} \,d\xi \\
     & \qquad + |x|^{-(\kappa+1)} \int_{\partial B(0,|x|^{-1})} \xii^{\eta-\kappa}\,dS(\xi)\\
     & \lesssim |x|^{-(n+\eta)}.
\end{split}
\]
Then:
\[\begin{split}
\|\mathscr{F}^{-1} (b(\xi)\,e^{-ta(\xi)})\|_{L^r}
     &\lesssim t^{-\frac{n+\eta}\theta}\,\Big(\int_{|x|\leq t^{\frac1\theta}}1dx\Big)^{\frac1r} + \Big(\int_{|x|\geq t^{\frac1\theta}} |x|^{-(n+\eta)r}\,dx\Big)^{\frac1r}\\
     & \lesssim t^{-\frac{n}\theta\left(1-\frac1r\right)-\frac\eta\theta},
\end{split} \]
provided that
\[ (n+\eta)r > n, \]
that is, \eqref{eq:reta} holds. If $\eta$ is an integer, then we just integrate by parts $\kappa=\eta+n-1$ times as in~\eqref{eq:intparts} and then twice as in~\eqref{eq:splitintparts} with $\delta=|x|^{-1}$.
\end{proof}
We then prove Lemma~\ref{lem:exp3}.
\begin{proof}[Proof of Lemma~\ref{lem:exp3}]
We proceed as in the proof of Lemmas~\ref{lem:exp} and~\ref{lem:exp2}, but now we notice that for any $|\gamma|\geq1$, we may estimate
\[ |\partial_\xi^\gamma (b(\xi)\,e^{-ta(\xi)})|\lesssim \big(\xii^{\eps-|\gamma|}+ t\,\xii^{\theta-|\gamma|}\big)\,e^{-\frac{c}2t\xii^\theta} \lesssim t^{\frac\eps\theta}\,\xii^{\eps-|\gamma|}\,,\]
due to $t\geq1$, where we assumed with no loss of generality that $\eps<1$ and $\eps\leq\theta$ (otherwise, it is sufficient to replace $\eps$ in the following by $\eps_1<\min\{1,\theta\}$). When we integrate by parts, we obtain
\[ |\mathscr{F}^{-1}(b(\xi)\,e^{-ta(\xi)})|\lesssim t^{\frac\eps\theta}\,|x|^{-(n+\eps)}, \]
so that
\[\begin{split}
\|\mathscr{F}^{-1} (b(\xi)\,e^{-ta(\xi)})\|_{L^r}
     &\lesssim t^{-\frac{n}\theta}\,\Big(\int_{|x|\leq t^{\frac1\theta}}1dx\Big)^{\frac1r} + t^{\frac\eps\theta}\,\Big(\int_{|x|\geq t^{\frac1\theta}}|x|^{-(n+\eps)r}\,dx\Big)^{\frac1r}\\
     & \lesssim t^{-\frac{n}\theta\left(1-\frac1r\right)},
\end{split} \]
as in the proof of Lemma~\ref{lem:exp}.
\end{proof}
Finally, we prove Lemma~\ref{lem:hexp3}.
\begin{proof}[Proof of Lemma~\ref{lem:hexp3}]
We first notice that for any $|\gamma|=n,n+1$, we have
\[ (ix)^\gamma \mathscr{F}^{-1}\phi = \mathscr{F}^{-1}(\partial_\xi^\gamma \phi)\in\mathcal C_0, \]
due to $\partial_\xi^\gamma \phi\in L^1$. On the one hand, letting $|\gamma|=n+1$, we find
\[ |\mathscr{F}^{-1}\phi(x)|\leq C\,|x|^{-(n+1)}. \]
On the other hand, letting $|\gamma|=n$ and $\gamma=\alpha+e_k$ for some $k$, where $|\alpha|=n-1$, we may integrate by parts
\[
\begin{split}
|(ix)^\gamma\mathscr{F}^{-1}\phi(x)|
     & \lesssim |x_k| \int_{B(0,|x|^{-1})} |\partial_\xi^\alpha \phi(\xi)|\, d\xi 
     + \int_{\R^n\setminus \bar B(0,|x|^{-1})} |\partial_\xi^\gamma \phi(\xi)| \,d\xi \\
     & \qquad + \int_{\partial B(0,|x|^{-1})} |\partial_\xi^\alpha \phi(\xi)|\,dS(\xi)\\
     & \lesssim |x|\,\int_{B(0,|x|^{-1})} \xii^{-\eps-(n-1)}\, d\xi 
     + \int_{\R^n\setminus \bar B(0,|x|^{-1})} \xii^{-\eps-n} \,d\xi \\
     & \qquad + \int_{\partial B(0,|x|^{-1})} \xii^{-\eps-(n-1)}\,dS(\xi)\\
     & \lesssim |x|^\eps.
\end{split}
\]
Therefore,
\[ |\mathscr{F}^{-1}\phi(x)|\leq C\,|x|^{-(n-\eps)}. \]
In turn, this proves that $\mathscr{F}^{-1}\phi\in L^1$. If~$\phi$ vanishes in a neighborhood of the origin, it is sufficient to notice that $\phi_j = i\xi_j\phi(\xi)/\xii$ still verify the assumptions of the first part of Lemma~\ref{lem:hexp3}, so that (see~\eqref{eq:H1norm} and~\eqref{eq:M1H1})
\[ \|\phi\|_{M(L^1,H^1)} = \|\mathscr{F}^{-1}\phi\|_{H^1} \leq \|\mathscr{F}^{-1}\phi\|_{L^1} + \sum_{j=1}^n \|\mathscr{F}^{-1}\phi_j\|_{L^1}<\infty. \]
\end{proof}

\section*{Declarations}

The authors declare that the data supporting the findings of this study are available within the paper.

The authors have no relevant financial or non-financial interests to disclose.

This paper has been partly realized during the stay of the second author at the Department of Mathematics of University of Bari in the period September-December 2022 and in the period June-August 2024, supported by the ``Visiting professor program'' of University of Bari. The first author is supported by PNR MUR Project CUP\_H93C22000450007, by ``PANDORA'' Project CUP\_H99J21017500006 and by INdAM-GNAMPA Project\\ CUP\_E55F22000270001. The second author is partially supported by ``Conselho Nacional de Desenvolvimento Cient\'ifico e Tecnol\'ogico (CNPq)" Grant 304408/2020-4.



\end{document}